\documentclass{amsart}
\usepackage[marginratio=1:1]{geometry}

\usepackage[hidelinks]{hyperref}
\usepackage{calc}
\newsavebox\CBox
\newcommand\hcancel[2][0.5pt]{%
  \ifmmode\sbox\CBox{$#2$}\else\sbox\CBox{#2}\fi%
  \makebox[0pt][l]{\usebox\CBox}%  
  \rule[0.5\ht\CBox-#1/2]{\wd\CBox}{#1}}
\usepackage{blindtext}
\usepackage{color}
\usepackage{hyperref,url}
\usepackage{float}
\usepackage{hyperref}
\usepackage{enumerate}
\usepackage{enumitem}   
\usepackage{bbm}
\usepackage{cancel}
\usepackage{mathrsfs}
\usepackage{colonequals}
\usepackage{pdfpages}

\usepackage{amsmath,amsfonts,amsthm,bm}
\usepackage{mathrsfs}  
\usepackage{xcolor}
\usepackage{amsfonts}
\usepackage{amssymb}
\usepackage{amsmath}
\usepackage{mathtools}
\usepackage{mathrsfs}  
\usepackage{bbm}
\usepackage{chngcntr}

\usepackage[normalem]{ulem}

\numberwithin{equation}{section}
\usepackage[toc,page]{appendix}
\usepackage{tikz-cd}
\theoremstyle{definition}

\theoremstyle{plain}
\newtheorem{teorema}{Theorem}[section]
\newtheorem{suposicao}{Hypothesis}

\newtheorem{lema}[teorema]{Lemma}
\newtheorem{proposicao}[teorema]{Proposition}

\newtheorem{step}{Step}
\counterwithin*{step}{teorema}

\theoremstyle{definition}
\newtheorem{definicao}{Definition}[section]
\newtheorem{notacao}[definicao]{Notation}

\theoremstyle{remark}

\newtheorem{remark}[teorema]{Remark}

\definecolor{roxo}{rgb}{0.44, 0.16, 0.39}
\definecolor{ao(english)}{rgb}{0.0, 0.5, 0.0}
\definecolor{dmagenta}{RGB}{139, 0, 139}
\definecolor{dgreen}{RGB}{0,90,0}
\definecolor{navy}{RGB}{0,0,128}

\usepackage{stackengine}

\def\R{\mathbb R}
\def\d{\mathrm d}
\def \d {\mathrm{d}}
\def \supp {\mathrm{supp}}

\def \cL {\mathcal{L}}

\definecolor{iblue}{RGB}{0, 35, 194}
\title[On the quasi-ergodicity of absorbing Markov chains]{On the quasi-ergodicity of absorbing Markov chains with unbounded transition densities, including random logistic maps with escape}
\author[M.M. Castro, V.P.H. Goverse, J.S.W. Lamb, M. Rasmussen]
{Matheus M. Castro$^*$,Vincent P.H. Goverse$^{*}$, Jeroen S.W. Lamb$^{*,\dagger,\ddagger}$, and Martin Rasmussen$^{*}$
}
\allowdisplaybreaks

\begin{document}



\maketitle
\vspace{-1.3cm}
\begin{center} \vspace{1cm}\small
   $^*$Department of Mathematics, Imperial College London, London, SW7 2AZ, United Kingdom.\\
   $^\dagger$International Research Center for Neurointelligence, The University of Tokyo, Tokyo, 113-0033, Japan.\\
  $^\ddagger$Centre for Applied Mathematics and Bioinformatics, Department of Mathematics and Natural Sciences, Gulf University for Science and Technology, Halwally, Kuwait.
\end{center}

\begin{abstract}

In this paper, we consider absorbing Markov chains $X_n$ admitting a quasi-stationary measure $\mu$ on $M$ where the transition kernel $\mathcal P$ admits an eigenfunction $0\leq \eta\in L^1(M,\mu)$. We find conditions on the transition densities of $\mathcal P$ with respect to $\mu$ which ensure that $\eta(x) \mu(\d x)$ is a quasi-ergodic measure for $X_n$ and that the Yaglom limit converges to the quasi-stationary measure $\mu$-almost surely.  We apply this result to the random logistic map $X_{n+1} = \omega_n X_n (1-X_n)$  absorbed at $\R \setminus [0,1],$ where $\omega_n$ is an i.i.d sequence of random variables uniformly distributed in $[a,b],$ for $1\leq a <4$ and $b>4.$

\end{abstract}

\section{Introduction}
Consider a family of transformations $\mathcal F = \{f_\omega :E\to E\}_{\omega \in \Delta}$, where $E$ is a metric space. Given a subset $M \subset E$ and endowing $\mathcal F$ with a probability measure, we aim to understand the statistical behaviour of the random dynamical system $$f^n (x;\omega_1, \ldots, \omega_n) := f_{\omega_n} \circ \ldots \circ f_{\omega_1}(x)$$ conditioned upon remaining in $M$.

Such a problem can be naturally modelled via the Markov chain $X_n = f_n(X_{n-1})$ with $f_n\in\mathcal F$ absorbed at $\partial := E\setminus M$, i.e. $X_n\in \partial$ implies $X_{n+1} \in \partial$. Statistical information for the above (conditioned) random dynamical system is then obtained by certain limiting distributions for the paths $X_n$. In the literature, such limiting distributions appear mainly in two forms. The first is the so-called \emph{Yaglom limit}
\begin{align}\label{Yag}
    \lim_{n\to \infty } \mathbb P[X_n \in A \mid X_0 = x, \tau >n],
\end{align}
where $x\in M$, $\tau =\min\{n\in\mathbb N\mid X_n\in \partial\}$ and $A$ is a measurable subset of $M$. The second one is the so-called \emph{quasi-ergodic limit}
\begin{align}\label{Qe}
   \lim_{n\to\infty}\mathbb E \left[\left.\frac{1}{n}\sum_{i=0}^{n-1} \mathbbm 1_A \circ X_n\, \right|\, X_0=x,\tau >n\right]. 
\end{align}

There are several contexts in which the Yaglom limit converges to a quasi-stationary measure. We recall that a probability measure $\mu$ on $M$ is called a \emph{quasi-stationary} measure for $X_n$ on $M$ if for every $n\in \mathbb N$, $\mu(\d x) := \mathbb P[X_n \in \d x\mid X_0 \sim \mu , \tau >n].$ On the other hand, the limit \eqref{Qe} is related to the existence of a quasi-ergodic measure for $X_n$ on $M$. A measure $\nu$ on $M$ is called a \emph{quasi-ergodic measure} for $X_n$ on $M$ if for every measurable subset $A$ of $M$ \eqref{Qe} converges to $\nu(A)$ for $\nu$-almost every $x\in M.$ 

In the literature \cite{benaim2021,QEM,Me,Champ,d,L2}, various sufficient conditions are presented for the existence and uniqueness of  quasi-stationary measures $\mu$ and quasi-ergodic measures $\nu$.  These conditions imply that $\nu \ll \mu$ and that the Radon–Nikodym derivative  ${\eta(x) =\nu(\d x)/\mu(\d x)}$ is an eigenfunction of $\mathcal P,$ where $\mathcal P$ is the transition kernel of $X_n$. The uniform convergence of the sequence $\left\{\mathcal P^n(\cdot,M)/\lambda^n\right\}_{n\in\mathbb N}$, where $\lambda := \int_M  \mathcal P(y,M)\mu(\d y)$,   to the eigenfunction $\eta$ plays a crucial role in the proofs.

{In this paper, we take a different approach. We set out to derive a quasi-ergodic measure starting from a quasi-stationary measure. 
The existence of quasi-stationary measures is a well-established problem (see \cite{survey} for a bibliography). 
Quasi-stationary measures arise as positive eigenmeasures of the operator $\mu \mapsto \int_M \mathcal P(x,\cdot )\mu(\d x)$  and an extensive literature exists on how to solve such eigenvalue problems \cite{MBanach, Krein, Schaefer}. Since quasi-ergodic measures do not admit such an approach, they are less well understood. Quasi-ergodic measures are important in the analysis of random dynamical systems, for instance, in the context of the recently established conditioned Lyapunov spectrum \cite{LEME, LExp}.

Inspired by these results, where quasi-ergodic measures can be expressed as a density over the quasi-stationary measures, we obtain natural conditions on the transition kernel $\mathcal P$ such that the existence of a quasi-ergodic measure becomes equivalent to solving an eigenvalue problem for $\mathcal P$ in $L^1.$ 
As a result, we considerably simplify the procedure of finding a quasi-ergodic measure. Furthermore, we also obtain, under aperiodicity conditions, that the Yaglom limit \eqref{Yag} converges to the quasi-stationary measure $\mu$-almost surely. }

As an application of our results, we characterise the limits \eqref{Yag} and \eqref{Qe} for the random logistic map $Y_{n+1} = \omega_n Y_n(1-Y_n)$ absorbed at $\R \setminus[0,1],$ with $\{\omega_n\}_{n\in\mathbb N}$ an i.i.d. sequence of random variables such that $\omega_0 \sim \mathrm{Unif}([a,b]),$ with $1\leq a <4$ and $b>4,$ where $\mathrm{Unif}([a,b])$ denotes the continuous uniform distribution in $[a,b]$.  The analysis of this system is challenging since its transition kernel presents a change of behaviour on the points $0$ and $1$. In particular, for every $x\in (0,1)$
$$\mathcal P(x,\d y) = \mathbb P[Y_n \in \d y \mid Y_0 =x] \ll \mathrm{Leb}(\d y),$$
while $\mathcal P(0,\d y ) = \mathcal P(1,\d y) = \delta_0(\d y).$ This implies that the transition densities $\mathcal P(x,\d y)$ explode when $x$ approaches the points $0$ and $1$. Consequently, the results in the literature \cite{benaim2021,QEM,Me,Champ,L2} cannot be applied, since $\mathcal P$ does not act as a compact operator on $\mathcal C^0(M)$ and $L^p(M)$, with $p\geq 1$. Hence, a more refined analysis is needed.

To overcome this issue, we consider $AM$-compact operators (see \cite[Appendix A]{GluckAsym}), a generalisation of compact operators. Inspired by the novel results on positive integral operators in \cite{gluckpostive,GluckAsym,JochenAper}, we analyse the action of $\mathcal P$ on $L^1(M,\mu)$. Since $\mathcal P$ is an integral operator, it is $AM$-compact, and we can establish its peripheral spectrum from which the asymptotic behaviour of $Y_n$ follows.

This paper is divided into six sections. In Section \ref{main}, the basic concepts of the theory of absorbing Markov chains are briefly recalled, and the main underlying hypotheses of this paper are defined (Hypotheses \ref{(K)} and \ref{(R)}), and the main results of this paper are stated (Theorems \ref{Thm:Logistic}, \ref{teorema2}, \ref{teorema3}  and \ref{teorema4}).  In Section \ref{4}, it is shown that  Hypothesis \ref{(K)} implies that $\mathcal P/\lambda$ is a mean ergodic operator. Section \ref{5} is dedicated to a brief presentation of Banach lattice theory, the definition of an $AM$-compact operator, and the proof of Theorem \ref{teorema2}.
In Section \ref{6}, we combine the results of the previous sections to prove Theorems \ref{teorema3} and \ref{teorema4}. Finally, in Section \ref{7}, we analyse the asymptotic behaviour of the random logistic map $Y_n$, introduced above, and prove Theorem \ref{Thm:Logistic}.

\section{Main results \label{main}}

Let $E$ be a metric space and $M$ a subspace of $E$. We aim to study Markov chains on $E$ conditioned upon remaining in the set $M$. With this objective in mind, we denote as $E_M$ the topological space $M \sqcup \partial$ generated by the topological basis
$$\mathcal T = \{A\cap M;\ A\ \text{is an open set of }E\}\sqcup \partial,$$
where $\sqcup$ denotes disjoint union.  In this paper, we assume that $$X:=\left(\Omega, \left\{\mathcal F_n\right\}_{n\in\mathbb N_0}, \left\{X_n\right\}_{n\in\mathbb N_0}, \left\{\mathcal P^n\right\}_{n\in\mathbb N_0}, \{\mathbb P_x\}_{x \in {{E_M}}}\right)$$ is a  Markov chain with state space ${E_M}$, in the sense of \cite[Definition III.1.1]{RW}, i.e., the pair $(\Omega,\{\mathcal F_n\}_{n\in\mathbb N})$ is a filtered space; $X_n$ is an $\mathcal F_n$-adapted process with state space ${E_M}$; $\mathcal P^n$ a time-homogeneous transition probability function of the process $X_n$ satisfying the usual measurability assumptions and Chapman-Kolmogorov equation; $\{\mathbb P_x\}_{x\in{E_M}}$ is a family of probability function satisfying $\mathbb P_x[X_0=x] = 1$ for every $x\in {E_M}$; and for all $m,n\in \mathbb N_0$, $x\in E_M$, and  every bounded measurable function $f$ on ${E_M}$, 
$$\mathbb E_x\left[f\circ X_{m+n}\mid \mathcal F_n\right] = ({\mathcal P}^m f)(X_n)\   \mathbb P_x\text{-almost surely}.$$

We assume that $X_n$ is a Markov chain that is absorbed at $\partial,$ meaning that ${\mathcal P}(\partial,\partial) = 1. $  In view of the above definitions, it is natural to define the stopping time
$$\tau(\omega) := \inf\{n\in\mathbb N; X_n(\omega) \not\in M\}.$$

 Throughout the paper, the following notations are used.

\begin{notacao}
 Given a probability measure $\mu$ on $M$ and  $p\in[1,\infty], $ we  denote ${L^p \left(M,{\mathscr B}(M),\mu\right)}$   as $L^p(M,\mu)$ and $\mathcal M(M)$ as the set of Borel signed-measures on $M.$ Moreover, we denote
${\mathbb P_\mu (\cdot) := \int_{M} \mathbb P_x (\cdot) \mu (\d x).}$

We denote as $\mathcal C^0(M) :=\{f:M\to \mathbb R; \ f \text{ is continuous}\}$ and $\mathcal F_b(M)$ as the set of bounded Borel measurable functions on $M$.  Given $f\in {\mathcal F_b}(M)$ write
$$\mathcal P^n f(x):= \mathcal P^n\left(\mathbbm 1_M   f \right)(x) = \int_M f(y) \mathcal P^n(x, \mathrm{d}y), $$
and by abuse of notation, denote
$$ f\circ X_n := \begin{cases} f\circ X_n,& \text{if } X_n\in M;\\ 0, & \text{if } X_n \notin M .\end{cases}$$

Given a sub $\sigma$-algebra $\mathscr F$ of $\mathscr B(M)$ and $f\in L^1(M,\mu),$ we denote $ \mathbb E_\mu [f\mid \mathscr F ] \in L^1(M,\mathcal F,\mu)$ as the conditioned conditional expectation of $f$ given $\mathscr F$, i.e. the unique function in $L^1(M,\mathscr F,\mu)$ such that
$$\int_F f(x) \mu(\d x) = \int_F \mathbb E_\mu [f \mid \mathscr F]\mu(\d x)\ \text{for every }F\in\mathscr F.$$

We define the sets
$$\mathcal M_+(M) = \{\mu\in\mathcal M(M);\mu(A) \geq 0,\ \text{for every }A\in\mathscr{B}(M)\}, $$
and
$$L^p_+ (M,\mu) = \{f\in L^p(M,\mu);\ f \geq 0 \ \mu\text{-almost surely}\},\ \text{for every }p\in[1,\infty].$$

Finally, given a Banach space $E,$ we say that the sequence $\{x_n\}_{n\in\mathbb N}\subset E$ converges in the weak topology to $x\in E$, if for every bounded linear functional $\phi \in E^*,$ $\lim_{n\to\infty } \phi(x_n) = \phi(x).$ Moreover, we say that the sequence $\{\phi_n\}_{n\in\mathbb N}\subset E^*$ converges in the weak$^*$ topology to $\phi$, if $\lim_{n\to\infty }\phi_n(x) = \phi(x)$ for every $x\in E.$

\end{notacao}

Since stationary measures do not capture the behaviour of $X_n$ before absorption, they become irrelevant when dealing with absorbing Markov chains. Due to this issue, it is necessary to extend the concept of stationary measures to quasi-stationary measures. Below, we recall the definition of a quasi-stationary measure. 

\begin{definicao}
A Borel measure $\mu$ on a metric space $M$ to be a \emph{quasi-stationary measure} for the Markov chain $X_n$ if
$$\mathbb P_\mu\left[X_n \in \cdot \mid \tau >n\right] =\mu(\cdot), \ {\   \text{for all}} \ n\in\mathbb N.$$
We call $\lambda = \int_M\mathcal P(x,M)\mu(\d x)$ {the} \emph{survival rate} of $\mu$.
\end{definicao}

Observe that if  $X_n$ admits a quasi-stationary measure $\mu$ on $M$ with survival rate $\lambda$, then $\mathcal P$ may be seen as a bounded linear operator in $  L^\infty (M,\mu)$. Moreover, since \begin{align}
    \int_M \mathcal P(x, A )\mu(\d x) = \lambda \mu(A)\ \text{for every }A\in\mathscr{B}(M),\label{l1ex}
\end{align}
and $L^\infty(M,\mu)$ is dense in $ L^1(M,\mu),$ the operator $\mathcal P$ can be naturally extended as a bounded linear operator in $L^1(M,\mu).$ 

While we have that ergodic stationary measures can be described in terms of Birkhoff averages for classical Markov chains. This is not true anymore when dealing with absorbing Markov chains,  meaning that quasi-stationary measures cannot be described in terms of conditioned Birkhoff averages. This obstruction motivates the definition of quasi-ergodic measures.

\begin{definicao}
A measure $\nu$ on $M$ is called a \emph{quasi-ergodic measure} if for every  $f\in {\mathcal F_b}(M),$
$$\lim_{n\to\infty}\mathbb E_x\left[\frac{1}{n} \sum_{i=0}^{n-1} f\circ X_i \hspace{0.1cm}\Bigg\vert \hspace{0.1cm} \tau >n\right] = \int_M f(y) \nu(\d y),\ {\   \mbox{for $\nu$-almost every}}\  x\in M. $$ \label{defqem}
\end{definicao}

% In this paper, we explore the intrinsic connections between quasi-stationary and quasi-ergodic measures for the absorbing Markov chain $X_n$ on $M$. 

% In this section we apply the above results to the 
One of the main objectives of this paper is to study the statistical asymptotic  behaviour of the  Markov chain $Y_{n+1}^{a,b} := \omega_n Y_{n}^{a,b} (1-Y_n^{a,b})$ absorbed at $\mathbb R\setminus [0,1]$, where $\{\omega_n\}_n$ is an i.i.d. sequence of random variables such that $\omega_0 \sim \mathrm{Unif}([a,b])$ with $1\leq a<4$ and $b>4$.

We mention that in the case in which $Y_n^{a,b}$ does not escape from the interval $[0,1]$, i.e.  $1\leq a < b\leq 4$, \cite[Theorem 2]{LogistMap} and \cite[Proposotion 9.5]{AleJan} shows that  $Y_n^{a,b}$ admits a unique stationary measure $\mu_{a,b}$ for $Y_n^{a,b}$ on $[0,1]$ such  $\mu_{a,b} ((0,1)) = 1.$ For dynamical considerations of random logistic maps and an analysis of the case were the sample space is finite, see \cite{RLM}.

The following theorem describes the asymptotic distribution of $Y_n^{a,b}$ conditioned upon survival when $1\leq a < 4<b$, also establishing the existence of quasi-stationary and quasi-ergodic measures for $Y_n^{a,b}$ on $[0,1]$.

\begin{teorema} \label{Thm:Logistic}
Consider $M=[0,1]$, $1\leq a<4<b$,  the Markov chain $Y_n^{a,b}$ on $\mathbb R_{M}$ absorbed at $\partial =\R\setminus M$ and  $\tau^{a,b}(\omega) = \min \{n\in\mathbb N; Y_n^{a,b} \in \mathbb R\setminus [0,1]\}.$
Then
\begin{enumerate}
    \item[${(i)}$] $Y_n^{a,b}$ admits a quasi-stationary measure $\mu_{a,b}$ with survival rate $\lambda_{a,b}$ such that $\supp(\mu_{a,b}) =[0,1]$ and $\mu_{a,b} \ll \mathrm{Leb}$, where  $\mathrm{Leb}$ denotes the Lebesgue measure on $[0,1].$
    \item[${(ii)}$] There exists $\eta_{a,b} \in L^1(M,\mu)$ such that $\mathcal P \eta_{a,b} = \lambda_{a,b} \eta_{a,b} $, $\|\eta_{a,b}\|_{L^1(M,\mu)} =1$ and $\eta_{a,b} >0$ $\mu_{a,b}$-almost surely.
    \item[${(iii)}$] For every $h\in L^\infty(M,\mathrm{Leb}),$   and $x\in (0,1)$ 
    $$ \lim_{n\to\infty} \mathbb E_x \left[\left. \frac{1}{n}\sum_{i=0}^{n-1} h\circ Y_i^{a,b}\, \right|\, \tau^{a,b} >n\right] = \int_M h(y) \eta_{a,b}(y) \mu_{a,b}(\d y).$$
    \item[${(iv)}$] For every $h\in L^\infty (M,\mathrm{Leb})$   and $x\in (0,1)$  
    $$\lim_{n\to \infty} \mathbb E_x\left[ h \circ Y_i^{a,b}\mid \tau^{a,b} >n\right]  = \int h(y) \mu_{a,b}(\d y). $$
\end{enumerate}

\end{teorema}
Theorem \ref{Thm:Logistic} is proved in Section \ref{logistic}.  Later, generalise the above result allowing values of $a$ in the interval $(0,1)$ (see Theorem \ref{Thm:ABLogistic}). However, this result relies on the technical assumption of $(a,b)$ being an admissible pair (see Definition \ref{AdmissiblePair}). We were not able to show the existence of quasi-stationary and quasi-ergodic measures for all values of $a\in [0,1),$ which can be seen from technical details in the inequalities of Proposition \ref{prop615} that is used in Step 2 of Lemma \ref{milagre3}. This technical obstruction is explained in Remark \ref{Problem}.

We use a more general setup for the proof of the above theorem. We present two incrementally restrictive hypotheses, Hypothesis  \ref{(K)} and \ref{(R)} which are satisfied by $Y^{a,b}_n$ for every $(a,b) \in [1,4)\times(4,\infty),$ that implies results similar to Theorem \ref{Thm:Logistic}  in different modes of convergence (see Theorems \ref{teorema2}, \ref{teorema3} and \ref{teorema4}).

In the following, we recall the definition of an integral operator.
\begin{definicao}
Let $p,q \in [1,\infty)$ and $(\Omega_1,\mathcal F_1,\mu_1),$ $(\Omega_2,\mathcal F_2,\mu_2)$ be measure spaces. We say that the bounded linear map $T:L^p(\Omega_1,\mathcal F_1,\mu_1)\to L^q(\Omega_2,\mathcal F_2,\mu_2)$ is an integral operator if there exists a measurable function $\kappa :\Omega_2\times \Omega_1\to \mathbb R$, called \emph{kernel function}, such that for every $f\in L^p(\Omega_1,\mathcal F_1,\mu_1)$ 
    $$\kappa (x,\cdot) f(\cdot) \in L^1(\Omega_1, \mu_1)\ \text{for }\mu_2\text{-almost every }x\in \Omega_2, $$
    and 
    \begin{align*}
          Tf(x)= \int_{\Omega_1} f(y) \kappa( x,y) \mu_1 (\d y),\mbox{ for $\mu_2$-almost every $x\in \Omega_2$.}
    \end{align*}
\end{definicao}

For a large class of Markov processes, it is common the existence of a probability function $\rho$ on $M$ such that 
$$\mathcal P(x, \d y) \ll \rho(\d y)\ \mbox{for $\rho$-almost every }x\in M. $$
In such systems, it is natural to seek quasi-stationary measures that are absolutely continuous with respect to $\rho$. In this situation and assuming that $\mu \ll \rho$ we have from \eqref{l1ex} that $\mathcal P: L^1(M,\mu)\to L^1(M,\mu)$ is an integral operator.

It is also natural to assume that the absorbing Markov chain $X_n$ is irreducible, i.e. if there exists $A \in \mathscr B(M)$ such that $\mu (\{\mathcal P(\cdot ,A) >0\}\triangle A) =0,$ then either $\mu(A) = 0$ or $\mu(X\setminus A)=0,$ where $\triangle$ denotes the symmetric difference of sets. In cases where $X_n$ is not irreducible, it is always possible to separate the state space into irreducible regions and analyse each region separately. 

The conditions discussed above are summarised in Hypothesis \ref{(K)}.

\begingroup
\renewcommand\thesuposicao{(H1)} 
\begin{suposicao} \label{(K)} Let $X_n$ be an absorbing Markov chain on $E_M$ absorbed at $\partial$. We say that $X_n$ fulfils Hypothesis \ref{(K)} if:

\begin{enumerate}
    \item[$\boldsymbol{\mathrm{(H1a)}}$]   There exists a quasi-stationary measure $\mu \in \mathcal M_+(M)$ for the Markov chain $X_n$ with survival rate $\lambda$.    
    
    \item[$\boldsymbol{\mathrm{(H1b)}}$]   There exists $\eta \in L^1_+(M,\mu)$ such that $\mathcal P \eta = \lambda \eta$ and $\|\eta\|_{L^1(M,\mu)}=1$.

    \item[$\boldsymbol{\mathrm{(H1c)}}$] The transition kernels $\mathcal P:L^1(M,\mu) \to L^1(M,\mu)$ is an integral operator with kernel function $\kappa:M\times M \to \mathbb R_+.$
    
    % There exists a measurable function $\kappa :M\times M\to \mathbb R_+$ such that for every $f\in L^1(M,\mu)$ 
    % $$\kappa (x,\cdot) f(\cdot) \in L^1(M,\mu)\ \text{for } \mu\text{-almost every} \ x\in M, $$
    % and 
    % \begin{align*}
    %      \mathcal Pf(x)= \int_M f(y) \kappa( x,y)\mu (\d y).
    % \end{align*}
    \item[$\boldsymbol{\mathrm{(H1d)}}$] for every $A \in \mathcal B(M)$ such that $0<\mu(A) < 1$ 
$$\int_{M\setminus A} \int_A \kappa (x,y) \mu(\d y) \mu(\d x) > 0, $$
i.e. if there exists $A\in \mathscr B(M)$ such that
$$\mu\left(\{\mathcal P\mathbbm 1_A >0\} \triangle A\right) = 0,$$
then either $\mu(M\setminus A) = 0$ or $\mu(A) = 0.$
\end{enumerate}

\end{suposicao}
\endgroup

We mention that given an absorbing Markov chain $X_n$ satisfying Hypothesis \ref{(K)}, we obtain from \cite[Lemma 4.2.9 and Example i) on Page 262]{MBanach} that $\eta(x) >0$ for $\mu$-almost every $x\in M.$

The theorem below implies that under Hypothesis \ref{(K)}, $\eta(x) \mu(\d x)$ is the only candidate for the quasi-ergodic measure for $X_n$ on $M.$ Moreover, it is also shown that such a hypothesis implies  the existence of a maximal $m\in\mathbb N,$ with the following properties:
\begin{itemize}
    \item there exists measurable sets $C_0,\ldots,C_{m-1} \subset M$ such that $M= C_0 \sqcup \ldots \sqcup C_{m-1};$
    \item for every $n\in\mathbb N$, $X_n \in C_{k\ (\mathrm{mod}\ m)},$ then $X_{n+1} \in C_{k+1\ (\mathrm{mod}\ m)}.$
\end{itemize}

% Under Hypothesis \ref{(K)} the theorem below shows that if $X_n$ admits a quasi-ergodic measure the only possibility is to be equal $\eta(x) \mu(\d x).$

\begin{teorema} \label{teorema2}

Let $X_n$ be an absorbing Markov chain fulfilling  Hypothesis \ref{(K)} then the following assertions hold:
\begin{enumerate}
    \item[(i)]There exist a natural number $m\in \mathbb N$ and sets $C_m:= C_0, C_1, \ldots, C_{m-1} \in \mathscr B(M)$ such that $\mu(C_i) = 1/m$ for every $i\in\{0,1\ldots, m-1\}$ and
$$\left\{\mathcal P \mathbbm 1_{C_i}  > 0\right\} \subset  C_{i-1}\ \text{for every }i\in \{0,1,\ldots,m-1\}.$$
    \item[(ii)] {For every $f\in L^1(M,\mu)$
    $$ \frac{1}{n} \sum_{i=0}^{n-1} \frac{1}{\lambda^i}\mathcal P^i f \xrightarrow{n\to\infty} \eta \int_M  f(y) \mu(\d y), $$
    in $L^1(M,\mu)$ and $\mu$-almost surely.}
 \item[(iii)]{There exist non-negative functions $g_0,g_1,\ldots,g_{m-1}\in L^1(M,\mu)$, satisfying $$\mathcal P g_{j} = \lambda g_{j-1\ \mathrm{(mod }\ m\mathrm)}\ \text{and}\ \|g_j\|_{L^1(M,\mu)}=1$$ 
 for every $j\in\{0,1,\ldots,n-1\}$, such that given $\ell\in \{0,1,\ldots,m-1\}$ and $h\in L^\infty(M,\mu)$ the following limit holds
 $$ \frac{1}{\lambda^{nm+\ell}} \mathcal P^{nm+\ell} h \xrightarrow{n\to\infty} \sum_{s=0}^{m-1}g_s \int_M h(x) \mu(\d x)$$
 in $L^1(M,\mu)$.}
\item[(iv)] If in addition, we assume that $M$ is a Polish space, then for every $h\in L^\infty(M,\mu)$ 

 \begin{align}
   \left(x\mapsto \mathbb E_x \left[\frac{1}{n} \sum_{i=0}^{n-1}h\circ X_i \mid \tau > n\right]\right) \xrightarrow{n\to\infty} \int_M h(y) \eta(y) \mu(\d y) \label{weakconv}
 \end{align}
in the $L^\infty(M,\mu)$-weak{${}^*$} topology (see \cite[Chapter 3.4]{Brezis} for the definition and the main properties of the weak$^{*}$ topology), in particular, we obtain that $\eqref{weakconv}$ also converges weakly in $L^1(M,\mu).$
\end{enumerate}
\end{teorema}

Theorem \ref{teorema2} is proved in Section \ref{a2}.

It is observed that Theorem  \ref{teorema2} $(iv)$ gives us $L^\infty(M,\mu)$-weak$^{*}$ convergence of \eqref{weakconv}, in order to guarantee such convergence in $L^\infty (M,\mu)$ we require an additional regularity  hypothesis (Hypothesis \ref{(R)}) on the kernel functions of the operator $\mathcal P$.

\begingroup
\renewcommand\thesuposicao{(H2)} 
\begin{suposicao}\label{(R)} Let $X_n$ be a Markov chain $X_n$ on $E_M$ absorbed at $\partial$. We say that $X_n$ fulfils Hypothesis \ref{(R)} if:
\begin{enumerate}
    \item $X_n$ fulfils Hypothesis \ref{(K)}; and
    \item for $\mu$-almost every point $x\in M$, $\kappa(x,\cdot)\in L^\infty(M,\mu)$.  Equivalently, since $\mu$ is an inner regular measure \cite[Proposition A.3.2.]{OV3},  there exists a sequence of nested compact sets $\{K_i\}_{i\in\mathbb N}$ such that $\mu\left(\bigcup_{i\in \mathbb N} K_i \right) = 1,$ and for every $i\in \mathbb N,$ 
     $$  \underset{(x,y) \in K_i\times M}{\mathrm{ess~sup}^{\mu\otimes \mu}}\|\mathbbm 1_{K_i}(x) \kappa (x,y)\| <\infty.$$
\end{enumerate}

\end{suposicao}
\endgroup

We mention that, in practice, once $\mathrm{(H1a)}$ and $\mathrm{(H1b)}$ are verified, then  $\mathrm{(H1c)}$, $\mathrm{(H1d)}$ and \ref{(R)} can be readily verified.  We exemplify this in Section \ref{randomlogistic}  considering the absorbing Markov chain $Y^{a,b}_n$ (see the proof Theorem \ref{Thm:ABLogistic}).

In addition to quasi-stationary measures, the so-called \emph{Yaglom limit}
$$ \lim_{n\to\infty} \frac{\mathcal P^n(x,A)}{\mathcal P^n(x,M)},\ \text{for }x\in M\ \text{and}\ A\in\mathscr{B}(M),$$
provides an alternative perspective on the asymptotic behaviour of the paths $X_n$ conditioned on survival. Observe that for the Yaglom limit to exist, it is necessary that  $M$ does not exhibit a cyclic decomposition under $X_n$, i.e., that $m=1$ on item $(i)$ of Theorem \ref{teorema2}. 

The following two results provide conditions that ensure the existence of a quasi-ergodic measure for $X_n$ on $M$ and the convergence of the Yaglom-limit.

\begin{teorema}\label{teorema3}
Let $X_n$ be an absorbing Markov chain fulfilling  Hypothesis \ref{(K)}. If any of the following items hold:
\begin{itemize}
    \item[$(a)$] There exists $K>0$ such that $\mu(\{K<\eta\}) =1$ almost surely.
    \item[$(b)$] There exists  $ g\in L^1(M,\mu)$ such that
     $$ \frac{1}{\lambda^n} \mathcal P^n(\cdot,M) \leq g \ \text{for every }n\in \mathbb N. $$
     \item[$(c)$] The absorbing Markov chain $X_n$ fulfils Hypothesis \ref{(R)}.
     \end{itemize}
     
Then for every $h\in L^\infty(M,\mu)$
\begin{align}
    \lim_{n\to\infty} \mathbb E_x\left[\left. \frac{1}{n}\sum_{i=0}^{n-1} h \circ X_i\, \right|\, \tau >n\right] = \int_M h(y) \eta(y) \mu(\d y), \mbox{ for $\mu$-almost every $x\in M$}. \label{acredita1}
\end{align}

If, additionally,  $m=1$ in Theorem \ref{teorema2} $(i)$  then  
\begin{align}
     \lim_{n\to\infty} \mathbb E_x\left[h\circ X_n \mid \tau>n\right] = \lim_{n\to\infty} \frac{\mathcal P^n h(x)}{\mathcal P^n(x,M)} = \int_M h(y) \mu(\d y), \mbox{ for $\mu$-almost every $x\in M$}. \label{acredita2}
\end{align}

% Finally, if $X_n$ fulfils Assumption \ref{(R)} and the nested sets $\{K_i\}_{i\in\mathbb N}$ given in Assumption \ref{(R)} are compact, then the conclusions above hold replacing \emph{$\mu$-almost every $x\in M$} to \emph{for every $x\in \bigcup K_i.$}

\end{teorema}

Theorem \ref{teorema3} is proved in Section \ref{a3}.

The following theorem is a refinement of the previous theorem, allowing us to characterise the set where the convergence of \eqref{acredita1} and \eqref{acredita2} hold.

{
\begin{teorema}\label{teorema4}
Let $X_n$ be an absorbing Markov chain fulfilling  Hypothesis \ref{(R)} and suppose that $\mathcal P f|_{K_i}\in\mathcal C^0(K_i)$ for every $f\in L^1(M,\mu)$ and $i\in\mathbb N$, where $\{K_i\}_{i\in\mathbb N}$ is the nested sequence of compact sets given by the second part of Hypothesis \ref{(R)}. Then, given $h\in L^\infty(M,\mu)$, \eqref{acredita1} holds for every $x\in \left(\bigcup_{i\in\mathbb N} K_i\right)\cap \{\eta >0\}.$

In case that $m=1$ in Theorem \ref{teorema2} $(i)$, \eqref{acredita2} holds for every $x\in  \left(\bigcup_{i\in\mathbb N} K_i\right)\cap\{\eta>0\}$.
\end{teorema}}

Theorem \ref{teorema4} is proved in Section \ref{a3}.
{
\begin{remark}
Notice that Theorems \ref{teorema3} and \ref{teorema4} also hold in a non-escape context. This means that if $X_n$ is a Markov chain on the metric space $M$ without absorption, satisfying the following properties:
\begin{itemize}
    \item $\mu$ is ergodic stationary measure for $X_n$ on $M;$ 
    \item the transition kernel $\mathcal P:L^1(M,\mu)\to L^1(M,\mu)$ is an integral operator; and
    \item $X_n$ is aperiodic, i.e.  $m=1$ in Theorem \ref{teorema2} (i),
\end{itemize}
then for every $h \in L^\infty(M,\mu),$ $\displaystyle\lim_{n\to\infty} \mathcal P^n h = \int h \,\d \mu,$ $\mu$-almost surely. In particular, from \cite[Theorem 1 (ii)]{WeakMixing} we obtain that $X_n$ is a weak-mixing Markov chain.
\end{remark}
}

\section{Mean-ergodic operators\label{4}}

For classical dynamical systems and Markov processes, mean-ergodic operators provide a vast array of tools and techniques for analysing their statistical properties  \cite[Chapters 7, 8 and 10]{OPAET}. This section shows that this is also true for absorbing Markov chains.

In the following, we recall the definition of a mean ergodic operator.

\begin{definicao}
Let $(E,\|\cdot\|)$ be a Banach space, we say that $T:E\to E$ is a \emph{mean-ergodic operator} if there exists a projection $P:E\to E$ such that
$$\lim_{n\to\infty}\left\|\frac{1}{n} \sum_{i=0}^{n-1} T^i x - Px\right\| =0,\ \text{for every }x\in E. $$

Let $M$ be a metric space, $\rho$ a Borel probability measure on $X,$ and  $T:L^1(M,\rho)\to L^1(M,\rho).$ We denote 
$$\mathcal I(T,\rho) = \sigma\left(A\in \mathscr B(M); T^*\mathbbm 1_A = \mathbbm 1_A\right), $$
where $T^*:L^\infty(M,\rho)\to L^\infty(M,\rho)$ is the dual operator of $T,$ i.e. the unique bounded automorphism on $L^\infty(M,\rho)$ such that 
$$\int_M (Tf)(x) h(x)~\rho(\d x) = \int f(x) (T^*h)(x)~\rho(\d x)\ \mbox{for every $f\in L^1(M,\rho)$ and $h\in L^\infty(M,\rho)$}. $$
\end{definicao}

Our results are highly dependent on the following two propositions.

\begin{proposicao}[{\cite[Theorem 3.3.5]{ET}} and {\cite[Corollary V.8.1]{Schaefer}}]\label{ET}
Let $M$ be a metric space and $\rho$ be a probability measure on $M$ and $T:L^1(M,\rho) \to L^1(M,\rho)$ be a linear operator such that $\|T\| = 1.$ Assume that there exists $\eta \in L^1(M,\rho)$ satisfying $T\eta = \eta$ and $\rho(\{\eta>0\}) = 1.$ Then,
\begin{enumerate}
    \item[(i)] For every $f\in L^1(M,\rho)$
    $${\lim_{n\to \infty} \frac{1}{n} \sum_{i=0}^{n-1}T ^i f = \eta \frac{\mathbb E_\rho[f \mid \mathcal I(T,\rho)]}{\mathbb E_\rho[ \eta \mid \mathcal I(T,\rho)]}\ \mbox{$\rho$-almost surely}}.$$
    \item[(ii)] The operator $T$ is mean-ergodic.
\end{enumerate}
\end{proposicao}\label{a1}

While Hypothesis \ref{(K)} does not imply that $\mathcal P$ is a compact operator, the proposition shows that given $f\in L^1(M,\mu)$ the orbit  $\left\{\frac{1}{\lambda^n}\mathcal P^n f\right\}_{n\in\mathbb N}$ is weakly precompact.

\begin{proposicao}\label{wcompact}
Suppose that the Markov process $X_n$ satisfies Hypothesis \ref{(K)}, then for every $f\in L^1(M,\mu)$ the sequence 
$$\left\{\frac{1}{\lambda^i}\mathcal P^i f\right\}_{i\in\mathbb N}\ \mbox{ is weakly-$L^1(M,\mu)$ precompact.}$$ 
\end{proposicao}
\begin{proof}

Let $f\in L^1_+(M,\mu)$. Note that for every $i,m\in\mathbb N,$
$$0\leq  \frac{1}{\lambda^i}\mathcal P^i f  \leq  \frac{1}{\lambda^i}\mathcal P^i (f- f\wedge m\eta) +  \frac{1}{\lambda^i}\mathcal P^i (f\wedge m\eta)\leq  \frac{1}{\lambda^i}\mathcal P^i (f- f\wedge m\eta) +   m\eta,$$
where given two function $f_1,f_2$, we define $f_1\wedge f_2:= \min\{f_1,f_2\}$. Since $\left\|\frac{1}{\lambda^i} \mathcal P^i f\right\|_{L^1(M,\mu)} = \|f\|_{L^1(M,\mu)}$ for every $i\in \mathbb N$ and  $ m \eta \wedge f \xrightarrow{m\to\infty} f$ in $L^1(M,\mu)$ and $\mu$-almost surely, we can obtain that for every $\varepsilon>0$ there exists $\delta>0$ such that
$$\frac{1}{\lambda^i}\int_A \mathcal P^i f(x) \mu(\d x) <\varepsilon\ \text{if }\mu(A)<\delta\ \text{for every }i\in\mathbb N.    $$
From \cite[Page 87, item 3]{CLM} we conclude that $\left\{\frac{1}{\lambda^n}\mathcal P^n f\right\}_{n\in\mathbb N}$ is  weakly-$L^1(M,\mu)$ precompact.
\end{proof}

\section{\texorpdfstring{$AM$}{AM}-Compact Operators}
\label{5}
Observe that under Hypothesis \ref{(K)}, the operator $ \frac{1}{\lambda}\mathcal P :L^1(M,\mu) \to L^1(M,\mu)$, is well behaved in a functional analytical point of view. Namely, $\frac{1}{\lambda}\mathcal P$ is a positive integral operator whose orbits are weakly compact. The theory of Banach lattices provides powerful tools for studying the spectrum of such operators. In the following two paragraphs, we recall the definition of a Banach lattice (we follow the definitions provided in \cite[Chapter 2]{MBanach} and \cite[Chapter 2]{Schaefer}).

Given $(L,\leq)$ a partially ordered set and a set $B\subset L,$ we define, if exists
$$\sup B = \min\{\ell\in L;\ b\leq \ell,{\   \text{for all}}\ b\in B\} $$
and
$$\inf B = \max\{\ell\in L;\ \ell \leq b ,{\   \text{for all}}\ b\in B\}. $$
With the above definitions, we say that $L$ is a \emph{lattice}, if for every $f_1,f_2\in L,$
$$f_1\lor f_2 := \sup \{f_1,f_2\} ,\ f_1\land f_2:= \inf \{f_1,f_2\}  $$
exists. Additionally,  in the case that $L$ is a vector space and the lattice $(L,\leq)$ satisfies
$$f_1\leq f_2 \Rightarrow f_1+f_3\leq f_2+f_3,\ {\   \text{for all}}\  f_3\in L,\ \text{and} $$
$$f_1\leq f_2 \Rightarrow \alpha f_1 \leq \alpha f_2,\ {\   \text{for all}\ } \alpha>0, $$
then $(L,\leq)$ is called \emph{vector lattice}. Finally, if  $(L,\|\cdot\|)$ is a Banach space and the vector lattice  $(L,\leq )$ satisfies
$$|f_1|\leq |f_2|\Rightarrow \|f_1\|\leq \|f_2\|,$$
where $|f_1| := f_1\lor (-f_1),$ then the triple $(L,\leq,\|\cdot\|)$ is called a 
\emph{Banach lattice}. When the context is clear, we denote the Banach lattice $(L,\leq,\|\cdot\|)$ simply by $L$.

In this paper, we use two fundamental notions from Banach lattice theory. The first one is that of an \emph{ideal} of a Banach lattice, and the second one is that of an \emph{irreducible operator} on a Banach lattice. A vector subspace $I \subset L$, is called an \emph{ideal} if, for every $f_1,f_2\in L$ such that $f_2\in I$ and $|f_1|\leq |f_2|,$ we have $f_1\in I$. Finally, a positive linear operator $T:L\to L$ is called \emph{irreducible} if,  $\{0\}$ and $L$ are the only $T$-invariant closed ideals of $T$.

\label{AM}
 The theory of $AM$-compact operators provides a generalisation to the theory of compact operators. $AM$-compact operators are considerably more general than compact operators and possess a sufficient degree of regularity. In the following, we recall the definition of an $AM$-compact operator.

\begin{definicao}
Let $E$ be a Banach lattice and $Y$ a Banach space. A linear operator $T:E \to Y$ is called $AM$-compact if for every $x_1,x_2\in E$, $T([x_1,x_2])$ is precompact in $Y$, where ${[x_1,x_2] :=\{y\in E;\ x_1\leq y \leq x_2\}.}$

\end{definicao}

The following result shows us that all positive integral operators are $AM$-compact.

\begin{teorema}[{\cite[Proposition A.5]{gluckpostive}}]
Let $(\Omega_1,\mu_1)$ and $(\Omega_2,\mu_2)$ the $\sigma$-finite measure spaces and $p,q\in [1,\infty)$. Let $T:L^p(\Omega_1,\mu_1) \to L^q (\Omega_2,\mu_2)$ be a positive bounded integral operator, then $T$ is an $AM$-compact operator. \label{amcompact}
\end{teorema}

 Given $f\in L^1(M,\mu)$, the key to our results is to understand the asymptotic behaviour of the sequence $\left\{\frac{1}{\lambda^n}\mathcal P^n f\right\}_{n\in \mathbb N}$. It turns out that the behaviour of this sequence has a strong connection with the peripheral spectrum of $\mathcal P$. In this way, we denote $L^1(M,\mu;\mathbb C) := L^1(M,\mu) \oplus i L^1(M,\mu)$ and linearly extend the operator $\mathcal P$ to the Banach space $L^1(M,\mu;\mathbb C)$.

Here, we summarise the spectral properties implied by Hypothesis \ref{(K)}.

\begin{proposicao}\label{e}
Let $X_n$ be an absorbing Markov chain satisfying Hypothesis \ref{(K)}. Then 
\begin{enumerate}
    \item[(i)]for every $f\in L^1(M,\mu),$
\begin{align*}
    {\frac{1}{n}\sum_{i=0}^{n-1}\frac{1}{\lambda^i}\mathcal P^i f \xrightarrow{n\to\infty} \eta \int_M f(y)\mu(\d y),}
\end{align*}
in $L^1(M,\mu)$ and $\mu$-almost surely.

\item[(ii)] There exists a decomposition  ${L^1(M,\mu; \mathbb C)=  E_{\mathrm{rev}} \oplus  E_{\mathrm{aws}}},$ such that $E_{\mathrm{rev}}$ and $E_{\mathrm{aws}}$ are $\mathcal P$-invariant,
$$E_{\mathrm{rev}}= \mathrm{span} \left\{f \in L^1(M,\mu;\mathbb C);\ \frac{1}{\lambda}\mathcal P f = e^{\frac{2\pi j i}{m} } f\ \text{for some }j\in\{0,1,\ldots,m-1\} \right\},$$
and 
$$E_{\mathrm{aws}}= \left\{f \in L^1(M,\mu; \mathbb C); \frac{1}{\lambda^i} \mathcal P^i f\xrightarrow{n\to \infty }0, \ \text{in }L^1(M,\mu)\right\}.$$

Moreover, $$\mathrm{dim}\ker\left(\frac{1}{\lambda}\mathcal P - e^{\frac{2 \pi j i}{m}}\mathrm{Id}\right) = 1 \ \text{for every }j\in\{0,1,\ldots,m-1\}.$$

\end{enumerate}

\end{proposicao}
\begin{proof}
$(i)$: From Proposition \ref{ET} it is enough to show that if $A \in \mathcal I(\mathcal P/\lambda, \mu)$ then either ${\mu(A) = 0 }$ or $\mu(A) = 1.$ In order to see this, Let $A\in \mathscr B(M)$ such that $\frac{1}{\lambda} \mathcal P^*\mathbbm 1_A = \mathbbm 1_A,$
then
\begin{align*}
    0 &= \mu( A \cap (M\setminus A) ) = \int_M  \mathbbm 1_A(x) \mathbbm 1_{M\setminus A}(x)\mu(\d x) =\frac{1}{\lambda }\int_M  \mathcal P^*\mathbbm 1_A(x) \mathbbm 1_{M\setminus A}(x)\mu(\d x)\\
    &= \frac{1}{\lambda }\int_M  \mathbbm 1_A(x) \mathcal P\mathbbm 1_{M\setminus A}(x)\mu(\d x) =\frac{1}{\lambda }\int_A  \int_{M\setminus A}\kappa(x,y) \mu( \d y) \mu(\d x).
\end{align*}
From \ref{(K)} we obtain that either $\mu(A) = 1$ or $\mu(A)=0$.

$(ii)$: From Propositions \ref{wcompact} and \ref{amcompact} we have the semigroup $\left\{\frac{1}{\lambda^n} \mathcal P^n \right\}_{n\in\mathbb N}$ fulfils the standard assumptions of \cite[Section 6]{GluckAsym}. Combining \cite[Propositions 4.3, Theorem 2.2]{GluckAsym} and \cite[Proposition 16.27 and Corollary 16.32]{OPAET} we obtain that
$$E_{\mathrm{rev}} = \overline{\left\{f \in L^1(M,\mu; C); \frac{1}{\lambda}\mathcal P f = e^{2 i \pi \theta} f, \ \text{for some }\theta\in \mathbb R\right\}},$$
and
$$E_{\mathrm{aws}} = \left\{f\in L^1(M,\mu; \mathbb C); \frac{1}{\lambda^i} \mathcal P^i f \to 0 \ \text{in }L^1(M,\mu)\right\}.$$

Applying \cite[Theorem 6.1 a)]{GluckAsym}, we obtain that if $   \lambda e^{2 i \pi \theta} \in \sigma_{\mathrm{pnt}}(\mathcal  P)$ then   $\theta \in \mathbb Q.$ Observe \ref{(K)} implies that $\mathcal P/\lambda$ is an irreducible operator \cite[Example i), Page 262]{MBanach}. From \cite[Theorem 6.1. b)]{GluckAsym} we obtain that $\mathcal P/\lambda$ has only finitely many unimodular eigenvalues. Finally, from \cite[Theorem 4.2.13 iii)]{MBanach} (taking $x' = 1$) the proof is finished.

\end{proof}

Let $\sigma_{\mathrm{pnt}}\left(\frac{1}{\lambda}\mathcal P\right) := \left\{ \widetilde \lambda \in \mathbb C;\ \mbox{ there exists $h\in L^1(M),$ $\frac{1}{\lambda}\mathcal Ph = \widetilde \lambda h$}\right\}$ be  the point spectrum of the operator $\frac{1}{\lambda}\mathcal P.$ In \cite{Me,Esperanca}, it is shown that  the cardinality of $\mathbb S^1 \cap \sigma_{\mathrm{pnt}}\left(\frac{1}{\lambda}\mathcal P\right)$ is intrinsically connected with the existence a possible periodic behaviour of $X_n$ in a suitable partition of $M$. This remains true under Hypothesis \ref{(K)}, and such periodic behaviour is established in Lemmas \ref{peakblinders} and \ref{ciclying1}.

\begin{definicao}
Given an absorbing Markov chain $X_n$ that satisfies \ref{(K)}. We define ${m(X_n) := \#\left( \mathbb S^1 \cap \sigma_{\mathrm{pnt}}\left(\frac{1}{\lambda}\mathcal P\right)\right)},$ which is finite from Proposition \ref{e}.
\end{definicao}

From now on, we denote $m(X_n)$ simply as $m$.

\begin{lema}\label{nerdola}
Let $X_n$ be a Markov chain satisfying  Hypothesis \ref{(K)}. Then there exist eigenfunctions $g_0,g_1,\ldots,g_{m-1} \in  L^1_+(M,\mu)$ of $\mathcal P^m,$ such that, $\|g_j\|_{L^1(M,\mu)} =1,$ for every $j\in\{0,1,\ldots,m-1\},$ and $\mathrm{span}_{\mathbb C} (\{g_i\}_{i=0}^{m-1}) = \mathrm{ker}( \mathcal P^m - \lambda^m \mathrm{Id}). $

Moreover, the eigenfunctions $g_0,$ $g_1,$ $\ldots$, $g_{m-1}$  can be chosen in a way such they have disjoint support, i.e., defining $C_i = \{g_i>0\}$, for all $i\in\{0,\ldots,m-1\}, $ then
$ \mu\left( C_i\cap C_j \right) = 0,\ {\   \text{for all}} \ i\neq j.$

Furthermore, the family of sets $\{C_i\}_{i=0}^{m-1}$ satisfies $$ \mu\left( M\setminus \left( C_0 \sqcup C_1 \sqcup \ldots \sqcup C_{m-1}\right)\right) = 0.$$\label{peakblinders}
\end{lema}
 \begin{proof}

The proof follows from similar arguments and computations laid out in \cite[Proposition 6.9]{Me} with the following two adaptations:
\begin{enumerate}
    \item the space $\mathcal C^0(M)$ is replaced by $L^1(M,\mu);$ and
    \item the set equalities are replaced by the relation $\sim.$ Namely, given $A,B \in \mathscr B(M)$ are said to be equivalent, i.e. $A\sim B$, if $\mu(A\triangle B) = 0,$ where $A\triangle B := (A\setminus B) \cup (B\setminus A).$ 
\end{enumerate}

 \end{proof}

The proof of the following lemma is analogous to the proof \cite[Lemma 6.15]{Me}. 
\begin{lema}
Let $\{g_i\}_{i=0}^{m-1}\subset L^1_+(M,\mu)$ as in Proposition \ref{peakblinders}. Then, re-labeling the functions $g_i$'s if necessary, we have that $\mathcal P g_i = \lambda g_{i-1\ (\mathrm{mod}\ m)}$ and ${\mathcal P^* \mathbbm 1_{C_i} = \lambda \mathbbm 1_{C_{i+1\ (\mathrm{mod}\ m)}}}.$ In particular, this implies that
$$\{\mathcal P(x,C_i)>0\} \subset C_{i-1\ (\mathrm{mod}\ m)}\ \text{for every}\  i\in\{0,1,\ldots,m-1\}. $$
\label{ciclying1}
\end{lema}

From now on, we fix $\{g_i\}_{i=0}^{m-1}$ and $\{C_i\}_{i=0}^{m-1}$ as in Lemma \ref{ciclying1}, and we denote $g_j =g_{{j\ (\mathrm{mod}\ m)}}$ and $C_j = C_{{j\ (\mathrm{mod}\ m)}}$ for every $j\in\mathbb N$.

The following two lemmas are the last ingredients needed for the proof of Theorem \ref{teorema2}.

\begin{lema}\label{forgottenland}
Suppose the absorbing Markov chain $X_n$ satisfies Hypothesis \ref{(K)}. Then {$h\in L^\infty(M,\mu)$ and $\ell \in \{0,1,\ldots,m-1\}$
\begin{align}
    \frac{1}{\lambda^{mn+\ell}} \mathcal P^{mn+\ell}h  \xrightarrow{n\to\infty} \sum_{s=0}^{m-1} g_{s}\int_{C_{s+\ell}} h\, \d\mu \ \text{in }L^1(M,\mu).\label{supervitality1}
\end{align}

and
\begin{align}
    \frac{1}{nm+\ell} \sum_{i=0}^{mn+\ell-1} \frac{1}{\lambda^i} \mathcal P^i\left(h \frac{1}{\lambda^{n-i}}\mathcal P^{n-i}\mathbbm 1_{M}\right) \xrightarrow{n\to\infty} \sum_{s=0}^{m-1} \mu(C_{\ell +s})g_{s}\int_M h \eta\, \d\mu \text{ in } L^1(M,\mu).\label{supervitality2}
\end{align}

}
\end{lema}

{
\begin{proof}

From Proposition \ref{e}. There exists $\alpha_0,\ldots,\alpha_{m-1}\in\mathbb C$ and $v\in E_{\mathrm{aws}}$ such that
\begin{eqnarray}
    h = \sum_{s=0}^{m-1} \alpha_s g_s + v.\label{hsum}
\end{eqnarray}

\begin{step}\label{step:app_v_null}
   We show that $v \in E_{\mathrm{aws}}$ if and only if $\int_{C_i} v \, \d \mu = 0$ for every $i\in\{0,1,\ldots,m-1\}.$
\end{step}
    
Suppose first that $v \in E_{\mathrm{aws}}$. We claim that $\mathbbm 1_{C_i} v  \in E_{\mathrm{aws}}$ for all $i  \in \{0,1,\ldots, m-1\}$. Indeed, if $\mathbbm 1_{C_i} v \not \in E_{\mathrm{aws}}$, then $ v = \alpha_i g_i + w + \sum_{j\neq i}\mathbbm 1_{C_j}v$ with $\alpha_i \neq 0$ and $w \in E_{\mathrm{aws}}$. Since, $\mu(C_i \cap C_j) = 0$ for all $j \neq i$ $(\mathrm{mod}\ m)$, we obtain that $v \not \in E_{\mathrm{aws}}$. It follows that $$\left|\int_{C_i} v\,\d\mu\right| =\left|\int_{M} \mathbbm 1_{C_i}v\,\d\mu\right| =\left|\int_{M}\frac{1}{\lambda^n}\mathcal P^n( \mathbbm 1_{C_i}v )\,\d \mu \right|\xrightarrow[]{n\to\infty}0.$$

Reciprocally, assume that $\int_{C_i} v \,\d \mu = 0$ for every $i\in\{0,1,\ldots,k-1\}.$ Write $v = \sum_{i=0}^{k-1}\alpha_i g_i + w$, with $w\in  E_{\mathrm{aws}}$. Since $\int g_i \d \mu = 1$, we have that
$$\alpha_i = \int_{C_i} \alpha_i g_i \,\d\mu = \int_{C_i} \left(\sum_{j=0}^{k-1}\alpha_j g_j + w\right) \,\d\mu = \int_{C_i} v\, \d \mu = 0.$$
We obtain that $\alpha_i = 0$ for every $i\in\{0,1,\ldots,k-1\},$ which implies $v\in  E_{\mathrm{aws}}$.

\begin{step}
We show \eqref{supervitality1}.
\end{step}
Integrating \eqref{hsum} with respect to $C_i$, from Step $1$, we obtain that 
\begin{eqnarray}
    h = \sum_{s=0}^{m-1} g_s \int_{C_s} h\, \d \mu + v.\label{eq:dec1m}
\end{eqnarray}

Therefore, $$\frac{1}{\lambda^{nm + \ell}} \mathcal P^{nm+\ell} h = \sum_{s=0}^{m-1} g_{s-\ell} \int_{C_s} h\, \d \mu + \frac{1}{\lambda^{nm + \ell}} \mathcal P^{nm+\ell} v \xrightarrow[]{n\to\infty} \sum_{s=0}^{m-1} g_s\int_{C_{s+\ell}} h\, \d \mu,$$
in $L^1(M,\mu).$

\begin{step}
We show \eqref{supervitality2}.
\end{step}

As before, there exists $w\in E_{\mathrm{aws}},$ such that $\mathbbm 1_M = \sum_{s=0}^{m-1}\mu(C_s)g_s +w$. From a direct computation, we obtain
\begin{align*}
   &\frac{1}{mn+\ell} \sum_{i=0}^{mn+\ell-1} \frac{1}{\lambda^i} \mathcal P^i\left(h(\cdot) \frac{1}{\lambda^{mn+\ell-i}}\mathcal P^{mn+\ell-i}(\cdot,M)\right)\\
  &= \frac{1}{mn+\ell} \sum_{i=0}^{mn+\ell-1} \frac{\mathcal P^i}{\lambda^i} \left(h \frac{\mathcal P^{mn+\ell-i}}{\lambda^{mn+\ell-i}} \left(\sum_{s=0}^{m-1}\mu(C_s)g_s + v\right)\right)\\
   &= \sum_{s=0}^{m-1}\mu(C_s) \frac{1}{mn+\ell} \sum_{i=0}^{mn+\ell-1} \frac{\mathcal P^i}{\lambda^i}\left(h g_{s-mn-\ell+i}\right) + \frac{1}{mn+\ell} \sum_{i=0}^{mn+\ell-1} \frac{\mathcal P^i}{\lambda^i} \left(h \frac{\mathcal P^{mn+\ell-i}}{\lambda^{mn+\ell-i}} v\right)\\
   &= \sum_{s=0}^{m-1}\mu(C_s) \frac{1}{mn+\ell} \sum_{i=0}^{mn+\ell-1} \frac{\mathcal P^i}{\lambda^i}\left(h g_{s-\ell+i}\right) + \frac{1}{mn+\ell} \sum_{i=0}^{mn+\ell-1} \frac{\mathcal P^i}{\lambda^i} \left(h \frac{\mathcal P^{mn+\ell-i}}{\lambda^{mn+\ell-i}} v\right).
\end{align*}

On the one hand, we have that
\begin{align*}
  \left\|\frac{1}{mn+\ell} \sum_{i=0}^{mn+\ell-1} \frac{\mathcal P^i}{\lambda^i} \left(h \frac{\mathcal P^{mn+\ell-i}}{\lambda^{mn+\ell-i}} v\right)\right\|_{L^1(M,\mu)} &\leq \frac{1}{mn+\ell} \sum_{i=0}^{mn+\ell-1} \left\|\frac{\mathcal P^i}{\lambda^i} \left(h \frac{\mathcal P^{mn+\ell-i} }{\lambda^{mn+\ell-i}}v\right)\right\|_{L^1(M,\mu)}\\
  &\leq\frac{1}{mn+\ell} \sum_{i=0}^{mn+\ell-1} \left\|\frac{\mathcal P^i}{\lambda^i} \left(\left|h \frac{\mathcal P^{mn+\ell-i} }{\lambda^{mn+\ell-i}}v\right|\right)\right\|_{L^1(M,\mu)}\\
  &\leq \frac{\|h\|_\infty }{mn+\ell} \sum_{i=0}^{mn+\ell-1} \left\| \frac{\mathcal P^{mn+\ell-i} }{\lambda^{mn+\ell-i}}v\right\|_{L^1(M,\mu)}\xrightarrow[]{n\to\infty}0.
\end{align*}

On the other hand, Step $2$ yields that 
\begin{align*}
\sum_{s=0}^{m-1} \frac{\mu(C_s)}{mn+\ell} \sum_{i=0}^{mn+\ell-1} \frac{\mathcal P^i}{\lambda^i}\left(h g_{s-\ell+i}\right)=&\sum_{s=0}^{m-1} \frac{\mu(C_s)}{mn+\ell} \sum_{j=0}^{n-1} \frac{\mathcal P^{mj}}{\lambda^{m j}}\left(\sum_{i =0}^{m-1}\frac{\mathcal P^i}{\lambda^i}\left(h g_{s-\ell+i}\right)\right)\\
&+\sum_{s=0}^{m-1} \frac{\mu(C_s)}{mn+\ell}  \frac{\mathcal P^{mn-1}}{\lambda^{mn -1}} \left(\sum_{i=0}^{\ell}\frac{\mathcal P^i}{\lambda^i}\left(h g_{s-\ell+i}\right)\right)\\
&\xrightarrow[]{n\to\infty}  \sum_{s=0}^{m-1} \frac{\mu(C_s)}{m} \sum_{k=0}^{m-1} g_k \int_{C_k}\sum_{i =0}^{m-1}\frac{\mathcal P^i}{\lambda^i}\left(h g_{s-\ell+i}\right)\, \d\mu\\
&=\sum_{s=0}^{m-1} \frac{\mu(C_s)}{m} \sum_{k=0}^{m-1} g_k \sum_{i=0}^{m-1}\int_{C_{k+i}} h g_{s-\ell+i}\, \d\mu\\
&= \sum_{k=0}^{m-1} \mu(C_{\ell+k} )g_k \int_{M} h \eta\, \d\mu,
\end{align*}
in $L^1(M,\mu)$.

Therefore, $$\frac{1}{mn+\ell} \sum_{i=0}^{mn+\ell-1} \frac{1}{\lambda^i} \mathcal P^i\left(h(\cdot) \frac{1}{\lambda^{mn+\ell-i}}\mathcal P^{mn+\ell-i}(\cdot,M)\right)\xrightarrow[n\to\infty]{} \sum_{k=0}^{m-1} \mu(C_{\ell+k} )g_k \int_{M} h \eta\, \d\mu,$$ in $L^1(M,\mu),$ which proves the lemma.

\end{proof}

}

{

}
Now, we prove Theorem \ref{teorema2}.

\label{a2}
\begin{proof}[Proof of Theorem \ref{teorema2}]
Items $(i),(ii)$ and $(iii)$ follows directly from, respectively,  Propositions \ref{ciclying1},  , \ref{e} (i) and Lemma \ref{forgottenland}.

In the following we prove $(iv).$ Given $h\in L^\infty(M,\mu)$ define
\begin{align*}
    g_n(x) &:= \mathbb E_x\left[\frac{1}{n}\sum_{i=0}^{n-1} h\circ X_i \mid \tau >n\right]= \frac{\lambda^n}{\mathcal P^n(x,M)} \frac{1}{n}\sum_{i=0}^{n-1} \frac{1}{\lambda^i} \mathcal P^i\left(h(\cdot) \frac{1}{\lambda^{n-i}}\mathcal P^{n-i}(\cdot,M)\right)(x).
\end{align*}

It is clear that $\|g_n\|_{L^\infty(M,\mu)} \leq \|h\|_{L^\infty(M,\mu)}$ for every $n\in \mathbb N$. Since $M$ is a Polish space, from the Banach–Alaoglu theorem, we obtain that the space $$B_{\|\cdot\|_{L^\infty(M,\mu)}} (0,\|h\|_{\infty}) :=\left\{ g\in L^\infty(M,\mu); \|g\|_{L^\infty(M,\mu)}\leq \|h\|_{L^\infty(M,\mu)}\right\},$$
is compact metric space when endowed with the $L^\infty(M,\mu)$-weak${}^*$ topology. {Given $\ell\in\{0,\ldots,m-1\},$ let $\{g_{mn_k+\ell}\}_{n\in\mathbb N}$} be a $L^\infty(M,\mu)$-weak${}^*$ convergent subsequence of $\{g_n\}_{n\in \mathbb N}$, and denote its limit as $g$.

We show that $g = \int_M h \eta\, \d \mu$ $\mu$-almost surely, which implies $(iv).$ Observe that given $A\in\mathscr B(M),$ from Lemma \ref{forgottenland} we obtain that
{
\begin{align*}
     \int_A g \left[\sum_{s=0}^{m-1}\mu(C_{s+\ell}) g_s\right] \d \mu &= \lim_{k\to\infty} \int_M g_{m n_k+\ell}\frac{\mathcal P^{m n_k+\ell}\mathbbm 1_M}{\lambda^{m n_k+\ell}}  \mathbbm 1_A \, \d\mu \\
    &=\lim_{k\to\infty} \int_M \frac{1}{mn_k +\ell}\sum_{i=0}^{m n_k +\ell -1} \frac{\mathcal P^i}{\lambda^i} \left(h(\cdot) \frac{\mathcal P^{mn_k+\ell-i}}{\lambda^{n-i}}(\cdot,M) \right)(x)\mathbbm 1_A(x)\mu(\d x)\\
    &= \int_A\left(\int_M h\eta\, \d \mu\right) \left[\sum_{s=0}^{m-1}\mu(C_{s+\ell}) g_s\right] \d \mu.
\end{align*}
Since $\mu(\{\sum_{s=0}^{m-1}\mu(C_{s+\ell}) g_s>0\})=1$ and $\ell\in\{0,1,\ldots,m-1\}$ is arbitrary it follows that $g = \int_M h(x)\eta(x) \mu(\d x)$ $\mu$-almost surely.
}
\end{proof}

\section{Almost-sure convergence \label{6}}

In this section we strengthen the  $L^\infty(M,\mu)$-weak${}^*$ convergence given in Theorem \ref{teorema2} to $L^\infty(M,\mu)$ convergence.

Note that for every $n\in\mathbb N,$ $x\in M$ and $A\in\mathscr B(M)$
$$ \mathbb E_x \left[\left.\frac{1}{n} \sum_{i=0}^{n-1} \mathbbm 1_{A} \circ X_i \,\right|\, \tau >n \right] = \frac{\lambda^n}{ \mathcal P^n(x,M)} \frac{1}{n}\sum_{i=0}^{n-1} \frac{1}{\lambda^i} \mathcal P^i\left(\mathbbm 1_A(\cdot) \frac{1}{\lambda^{n-i}}\mathcal P^{n-i}(\cdot,M)\right)(x).$$
Therefore, to prove Theorem \ref{teorema3}, it is enough to find conditions that \eqref{supervitality1} and \eqref{supervitality2} converge almost surely. 
 
To prove Theorem \ref{teorema3} we need the following three propositions below. 

\begin{proposicao}[{\cite[Proposition 3.3.3.]{MBanach}}]\label{escolhi} Let $T:L^1(M,\mu)\to L^1(M,\mu)$ be a positive bounded integral operator. Then if $\{f_n\}_{n\in\mathbb N}\subset L^1(M,\mu)$  is a $L^1(M,\mu)$-order bounded sequence satisfying $f_n\to 0$ in $\mu$-measure as $n\to\infty$, then $Tf_n \to 0$ as $n\to \infty$ $\mu$-almost everywhere.
\end{proposicao}

\begin{proposicao}\label{rima1}
Let $X_n$ be an absorbing Markov chain satisfying Hypothesis \ref{(K)}. Suppose that one of the following items holds:
\begin{itemize}
    \item[$(a)$] There exists $K>0$ such that $\mu(\{K<\eta\}) =1$ almost surely.
    \item[$(b)$] There exists  $ g\in L^1(M,\mu)$ such that
     $$ \frac{1}{\lambda^n} \mathcal P^n(x,M) \leq g \ \text{for every }n\in \mathbb N. $$
     \item[$(c)$] The absorbing Markov chain $X_n$ fulfils Hypothesis \ref{(R)}.
     \end{itemize}
\end{proposicao}

{
Then, for every $h\in L^\infty(M,\mu)$ and $\ell\in\{0,1,\ldots,m-1\}$
\begin{align}
    \lim_{n\to \infty} \frac{1}{n}\sum_{i=0}^{mn+\ell-1}\frac{\mathcal P^i}{\lambda^i} \left(h\frac{\mathcal P^{mn+\ell-i}(\cdot,M)}{\lambda^{mn+\ell-i}} \right)  \xrightarrow{n\to\infty}  \sum_{s=0}^{m-1}\mu(C_{s+\ell})g_s\int_M h(y)\eta(y)\mu(\d y)\ \mu\text{-a.s.}.\label{elon1}
\end{align}

In addition,
\begin{align}
    \frac{1}{\lambda^{mn+\ell}} \mathcal P^{nm+\ell} h \xrightarrow{n\to\infty} \sum_{s=0}^{m-1}\mu(C_{s+\ell})g_s \int_M h(x) \mu(\d x)\, \mu\text{-a.s.}\label{elon2}.
\end{align} }
\begin{proof}
Observe that $(a)$ is a particular case of $(b)$. In fact, note that for every $x\in M,${
$\frac{1}{\lambda^n}\mathcal P^n(x,M)\leq \frac{1}{K}\frac{1}{\lambda^n}\mathcal P^n\eta(x) = \frac{\eta(x)}{K}$} which correspond to $(b)$
when setting $g := \eta/K.$ Now, we assume $(b)$. From Lemmas \ref{supervitality2} and \ref{forgottenland}, we obtain that \eqref{elon1} and \eqref{elon2} converge in probability. {Moreover, $(b)$ implies that for every $n\in\mathbb N,$ and for $\mu$-almost every $x \in M$
$$-\|h\|_{L^\infty(M,\mu)}g(x)\leq \frac{1}{\lambda^n} \mathcal P^{nm+\ell} h(x)\leq \|h\|_{L^\infty(M,\mu)}g(x)  $$
and
$$-\|h\|_{L^\infty(M,\mu)} g(x) \leq \frac{1}{nm+\ell}\sum_{i=0}^{nm+\ell-1}\frac{1}{\lambda^i} \mathcal P^i\left(h(\cdot)\frac{1}{\lambda^{{nm+\ell}-i}} \mathcal P^{{nm+\ell}-i}(\cdot,M)\right)(x) \leq \|h\|_{L^\infty(M,\mu)} g(x).$$
Therefore, Proposition $\ref{escolhi}$ implies the result.}

Now, we assume $(c)$. {Given $j\in\mathbb N$ consider the set
$$K_j = \{x\in M; \|k(x,\cdot)\|_{L^\infty(M,\mu)}\leq j\}. $$
It is clear that the map
\begin{align*}
    \mathcal G_j: L^1(M,\mu)&\to L^\infty(K_j,\mu)\\
    f&\mapsto  \frac{1}{\lambda} \mathbbm 1_{K_j}\mathcal P(f),
\end{align*}
is a bounded linear operator and $\|\mathcal G_j\|\leq j$. Given $h\in L^\infty(M,\mu)$, we obtain from Lemma \ref{supervitality2} that
\begin{align*}
    \mathbbm1_{K_j}\frac{1}{\lambda^{nm+\ell}} \mathcal P^{nm+\ell} h =  \mathcal G_j\left(\frac{1}{\lambda^{mn+\ell -1}}\mathcal P^{mn +\ell -1}\right)\xrightarrow[L^\infty(M,\mu)]{n\to\infty}& \sum_{s=0}^{m-1}\mu(C_{s+\ell-1}) \mathcal G_j (g_s) \int_M h \,\d\mu\\
    &=\mathbbm 1_{K_j}\sum_{s=0}^{m-1}\mu(C_{s+\ell})g_s \int_M h\, \d \mu.
\end{align*}

Similarly, from Lemma \ref{forgottenland}
\begin{align*}
   &\mathbbm 1_{K_j}  \frac{1}{nm+\ell-1}\sum_{i=0}^{nm+\ell-1}\frac{\mathcal P^i}{\lambda^i} \left(h\frac{\mathcal P^{nm+\ell-i}(\cdot,M)}{\lambda^{nm+\ell-i}}\right) \\
    &= \mathbbm 1_{K_j}  \frac{1}{nm+\ell-1}\sum_{i=1}^{nm+\ell-1}\frac{\mathcal P^i}{\lambda^i} \left(h\frac{\mathcal P^{nm+\ell-i}(\cdot,M)}{\lambda^{nm+\ell-i}}\right) +\mathbbm{1}_{K_j}\frac{1}{nm+\ell-1}h\frac{\mathcal P^{nm+\ell}}{\lambda^{mm+\ell}}(\cdot ,M)\\
    &= \mathcal G_j\left( \frac{1}{nm+\ell-1}\sum_{i=0}^{nm+\ell-2}\frac{\mathcal P^{i}}{\lambda^{i}} \left(h\frac{\mathcal P^{nm+\ell-1-i}(\cdot,M)}{\lambda^{nm+\ell-1-i}}\right) \right)+\frac{1}{nm+\ell-1}h \mathcal G_j\left(\frac{\mathcal P^{nm+\ell-1}}{\lambda^{mm+\ell-1}}(\cdot ,M) \right)\\
   &\xrightarrow[L^\infty(M,\mu)] {n\to\infty} \sum_{s=0}^{m-1}\mu(C_{s+\ell-1}) \mathcal G_j (g_s) \int_M h \eta \,\d\mu=\mathbbm 1_{K_j}\sum_{s=0}^{m-1}\mu(C_{s+\ell})g_s \int_M h \eta \,\d\mu.
\end{align*}

Since Hypothesis \ref{(R)} implies that $\mu(\bigcup_{j\geq 1} K_j) = 1$ we obtain the result.}
\end{proof}

\label{a3}

Finally, we prove Theorems \ref{teorema3} and \ref{teorema4}.
\begin{proof}[Proof of Theorem \ref{teorema3}]
{Given $h\in L^\infty(M,\mu)$, for every $n\in\mathbb N$ and $x\in M$ we obtain,
\begin{align}
   \mathbb E_x\left[\frac{1}{n} \sum_{i=0}^{n-1}h\circ X_i \mid \tau > n\right] = \frac{\displaystyle \frac{1}{n}\sum_{i=0}^{n-1} \frac{1}{\lambda^i}\mathcal P^i\left(h(\cdot) \frac{1}{\lambda^{n-i}}\mathcal P^{n-i}(\cdot,M) \right)(x)}{\displaystyle \frac{1}{\lambda^n} \mathcal P^n(x,M)}.\label{canyouhearme} 
\end{align} 

From Proposition \ref{rima1} we obtain that for $\mu$-almost every $x\in M$
$$\lim_{n\to\infty}\mathbb E_x\left[\frac{1}{n} \sum_{i=0}^{n-1}h\circ X_i \mid \tau > n\right]  =  \int_M h(x) \eta(x)\mu(\d x),$$
which proves the first part of the theorem.

In the case $m=1$ in Theorem \ref{teorema2}, Proposition \ref{rima1} implies that for $\mu$-almost every $x\in M$ 
$$\lim_{n\to\infty} \frac{\mathcal P^n h (x)}{ \mathcal P^n(x,M)} = \int_M h(y) \mu(\d y).$$}
\end{proof}

\begin{proof}[Proof of Theorem \ref{teorema4} ]
{Consider $h\in L^\infty(M,\mu)$, and let $j\in \mathbb N$. Observe that the operator
\begin{align*}
    \mathcal G_j:L^1 (M,\mu)&\to \mathcal C^0(K_j)\\
    g&\mapsto  \mathbbm{1}_{K_i}\frac{1}{\lambda}\mathcal Pg(x). 
\end{align*}
is bounded since $L^1(M,\mu)$ and $\mathcal C^0(K_j)$ are Banach lattices and $\mathcal G_j$ is a positive operator \cite[Theorem 5.3]{Schaefer}. Repeating the the proof of Proposition \ref{rima1}, we obtain that for every $\ell \in \{0,1,\ldots,m-1\},$
\begin{align}
   \frac{1}{nm+\ell}\sum_{i=0}^{nm+\ell-1}\frac{\mathcal P^i}{\lambda^i} \left(h(\cdot)\frac{\mathcal P^{mn+\ell-i}(\cdot,M)}{\lambda^{nm+\ell-i}} \right) 
    \xrightarrow[\mathcal C^0(K_j)]{n\to\infty}  \sum_{s=0}^{m-1}\mu(C_{s+\ell})g_s\int_M h(y)\eta(y)\mu(\d y).\label{everytime1}
\end{align}
and
\begin{align}
    \lim_{n\to\infty}\frac{1}{\lambda^{nm+\ell}}\mathcal P^{nm+\ell}(\cdot,M)  \xrightarrow[\mathcal C^0(K_j)]{n\to\infty}  \sum_{s=0}^{m-1}\mu(C_{s+\ell})g_s.\label{everytime2}
\end{align}

Since $\ell\in\{0,1,\ldots,m-1\}$ is arbitrary, equations \eqref{canyouhearme}, \eqref{everytime1} and \eqref{everytime2} yield that for every $x\in \bigcup_{i\in\mathbb N} K_i\cap \{\eta>0\}$ 
 $$ \lim_{n\to\infty }\mathbb E_x\left[\frac{1}{n} \sum_{i=0}^{n-1}h\circ X_i \mid \tau > n\right] = \int_M h(y) \eta (y)\mu(\d y).$$

Note that if $m=1$ in Theorem \ref{teorema2}. We conclude that
$$ \lim_{n\to\infty} \frac{\mathcal P^n h(x)}{\mathcal P^n(x,M)} = \int_M h(y) \mu(\d y),\ \text{for every }x\in\bigcup_{i\in\mathbb N} K_i \cap\{\eta>0\}.$$}

\end{proof}

\section{Random logistic map with escape\label{7}}
\label{randomlogistic}

In this section we analyse the Markov chain $Y_{n+1}^{a,b} = \omega_n Y_{n}^{a,b} (1-Y_n^{a,b})$ absorbed at $\partial = \mathbb R\setminus M,$ where $\{\omega_n\}_{n\in\mathbb N}$ is an i.i.d. sequence of random variables such that $\omega_n\sim\mathrm{Unif}([a,b])$, where $0<a<4<b$ and $M=[0,1].$ As before, for every $A\in\mathscr B(M)$ and $x\in M$ we denote
$$\mathcal P(x,A):= \mathbb P[ Y_1^{a,b} \in A \mid Y_0^{a,b} = x].$$

Clearly $\delta_0$ is a stationary measure for $Y_n^{a,b}$ on $[0,1]$. In the following, we provide conditions to show that $Y_n^{a,b}$ admits a non-trivial quasi-stationary measure on $[0,1]$, which we define as a quasi-stationary measure for $Y_n^{a,b}$ different from $\delta_0.$ For the sake of simplicity and in the interest of readability, we denote $Y^{a,b}_n$ simply as $Y_n$. Similarly, when the context is clear, we omit the $a,b$ superscript from future objects that depend on $a$ and $b$.

In the proposition below, we explicitly compute the transition functions of $Y_n.$

\begin{proposicao}\label{uk1}
Let $0\leq a \leq b $, and consider the absorbing Markov chain $Y_n^{a,b}$. Moreover, given  $f\in L^1(M,\mathrm{Leb}),$ $$\mathcal P f(x) = \frac{1}{(b-a)x(1-x)}\int_{ax(1-x)}^{bx(1-x)\wedge 1} f (y) \d y. $$
In the case that $f\in \mathcal C^0(M),$ then $\mathcal P f \in \mathcal C^0(M)$ and  $\mathcal P f (0) = \mathcal Pf(1) = f(0). $
\end{proposicao} 
\begin{proof}
Let $f\in L^1(M,\mathrm{Leb})$ by a direct computation,
\begin{align*}
    \mathcal P f(x) &= \frac{1}{ b-a} \int_a^b \mathbbm 1_{[0,1]} (\omega x(1-x)) f(\omega x (1-x)) \d \omega\\
    &=\frac{1}{(b-a)x(1-x)}\int_{ax(1-x)}^{bx(1-x)} \mathbbm 1_{[0,1]} (y) f(y) \d y\\
    &=\frac{1}{(b-a)x(1-x)}\int_{ax(1-x)}^{bx(1-x)\wedge 1} f(y) \d y .
\end{align*}
Now, consider $f\in \mathcal C^0(M).$ The above equation implies that $\mathcal Pf$ is continuous in $(0,1)$. For every $x\in (0,1)$, let us define the interval $J_x := [ax(1-x), bx(1-x)\land 1].$ It follows that for every $x\in (0,1/b)$, $\min_{y\in J_x} f(y) \leq \mathcal Pf(x) \leq \max_{y\in J_x} f(y).$ From the continuity of $f$ we obtain that $\lim_{x\to 0} \mathcal Pf = f(0).$ Since $\mathcal P(x) = \mathcal P (1-x)$ for every $x\in (0,1/2)$, if follows that $\lim_{x\to 1}\mathcal P(x) = f(0),$ implying that $\mathcal P f\in \mathcal C^0(M)$.

\end{proof}

The first step to apply Theorem \ref{teorema4} to $Y_n$ on $[0,1],$  is to show that $Y_n$ admits a quasi-stationary measure different from $\delta_0$ on $[0,1]$. 

Consider a measure $\mu \in \mathcal M (M),$ such that $\mu \ll \mathrm{Leb}(\mathrm{d}x)$  and define $ g:=\mu(\d x)/ \mathrm{Leb}(\d x).$ Note that
\begin{align*}
    \mathcal P^*(\mu)(A) &= \int_M \mathcal P(x,A) g(x) \d x \\
    &= \int_M \frac{1}{(b-a)x(1-x)}\int_0^{1}\mathbbm 1_{[ax(1-x),bx(1-x)]}(y) \mathbbm 1_A(y) g(x) \d y\d x\\
    &=  \int_{A}\int_{M}  \mathbbm 1_{[ax(1-x), b x(1-x)]} (y) \frac{  g(x)}{(b-a)x(1-x)} \, \d x \d y\\
    &=  \int_{A}\left(\int_{\alpha_-(y)}^{\alpha_+(y)} \frac{  g(x)}{(b-a)x(1-x)} \, \d x  - \int_{\beta_-(y\wedge a/4)}^{\beta_+(y\wedge a/4)} \frac{  g(x)}{(b-a)x(1-x)} \, \d x \right)\d y,
\end{align*}
where
$$\alpha_{\pm}(x) := \frac{1}{2} \pm \frac{1}{2}\sqrt{1 - \frac{4}{b}x}\ \text{and } \beta_{\pm}(x) := \frac{1}{2} \pm \frac{1}{2}\sqrt{1 - \frac{4}{a}x}.$$

The above observation motivates the definition of the stochastic transfer operator,
\begin{align}
\mathcal L: L^1([0,1],\mathrm{Leb}) &\to L^1([0,1],\mathrm{Leb}) \label{stransfer}\\   
g&\mapsto \left(x\mapsto  \int_{\alpha_-(x)}^{\alpha_+(x)} \frac{  g(y)}{(b-a)y(1-y)} \, \d y  - \int_{\beta_-(x\wedge a/4)}^{\beta_+(x\wedge a/4)} \frac{  g(y)}{(b-a)y(1-y)} \, \d y \right), \nonumber
\end{align}
note that  $\mathcal L$ is a well defined linear operator since for every $g\in L^1([0,1],\mathrm{Leb})$
$$ \left\|\mathcal L (g)\right\|_{L^1(M,\mathrm{Leb})} := \int_0^1 \left|\mathcal L g(x)\right| \d x \leq \int_M \mathcal P(x,M) |g(y)| \d y \leq \|g\|_{L^1(M,\mathrm{Leb})}. $$

The following two propositions summarise the above comments and show that $\mathcal L$ is well defined as an automorphism in $L^p(M,\mathrm{Leb})$ for every $p\in[1,\infty].$ For the following result see for instance  \cite[Section 5]{AleJan}.

\begin{proposicao}
A probability measure $\mu \in \mathcal M_+(M)\setminus\{\delta_0\}$ on $[0,1]$ is a quasi-stationary measure for $Y_n$ if, and only if, $\mu(\d x) \ll \mathrm{Leb}(\d x)$ and there exists $0<\lambda<1,$ such that
$$ \mathcal L \frac{\mu(\d x) }{\mathrm{Leb}(\d x)} = \lambda \frac{\mu(\d x) }{\mathrm{Leb}(\d x)}. $$
\end{proposicao}

% \begin{proof}
% Using Proposition \ref{uk1}, it is readily verified if $\mu$ is a quasi-stationary measure for $Y_n$ on $[0,1]$ different of $\delta_0$ then $\mu \ll \mathrm{Leb}$ and
% $\lambda = \int_M \mathcal P(x,M)\mu <1.$ Moreover, the computations above $\eqref{stransfer}$ show that $\mu$ is a quasi stationary measure if, and only if
% $$ \mathcal L \frac{\mu(\d x) }{\mathrm{Leb}(\d x)} = \lambda \frac{\mu(\d x) }{\mathrm{Leb}(\d x)},$$
% for some $\lambda<1.$
% \end{proof}

\begin{proposicao}
For every $p\in [1,\infty]$ the operators 
$$\left.\mathcal P\right|_{L^p([0,1],\mathrm{Leb})}, \left.\mathcal L\right|_{L^p([0,1],\mathrm{Leb})}: L^p([0,1]) \to  L^p([0,1]).  $$
are well defined and bounded.
\end{proposicao}

\begin{proof}
By a direct computation, one can check that
$$ \mathcal L \mathbbm 1_M (x) = 
\begin{cases}
\displaystyle\frac{4}{b-a}\left(\displaystyle\tanh^{-1}\left(\sqrt{1-\frac{4 x}{b}}\right)-\tanh ^{-1}\left(\sqrt{1-\frac{4x}{a} }\right)\right)&, \displaystyle\text{if }0\leq x\leq \frac{a}{4}\\
\displaystyle\frac{4}{b-a}\tanh^{-1}\left(\sqrt{1-\frac{4 x}{b}}\right)&,\displaystyle \text{if }\frac{a}{4}\leq x\leq 1.\\
\end{cases}$$
implying that
$$\|\mathcal L\|_{L^\infty(M,\mathrm{Leb})} = \frac{4}{b-a} \tanh^{-1}\left(\sqrt{1-\frac{a}{b}}\right).$$

Since $\|\mathcal L\|_{L^1(M,\mathrm{Leb})} \leq 1,$ by the Riesz–Thorin interpolation theorem \cite[Theorem 6.27]{Folland}
$$\|\mathcal L\|_{L^p(M,\mathrm{Leb})} < \infty, \ \text{for all } p\in \left[ 1 , \infty\right]. $$

For the operator $\mathcal P,$ note that for every $0\leq f \in L^1 ([0,1])$ ($1\leq p \leq \infty$),
$$\int_0^1 \mathcal P f(x)\, \d x = \int_0^1 f(x) \mathcal L \mathbbm 1_M(x) \, \d x \leq \|\mathcal L\|_{L^\infty(M,\mathrm{Leb})} \|f\|_{L^1(M,\mathrm{Leb})},$$
showing that $\|\mathcal P\|_{L^1(M,\mathrm{Leb})} \leq \|\mathcal L\|_{L^\infty(M,\mathrm{Leb})}<\infty.$ Using that $\|\mathcal P\|_{L^\infty(M,\mathrm{Leb})} \leq 1$ we have again by the Riesz-Thorin interpolation theorem that $\|\mathcal P\|_{L^p(M,\mathrm{Leb})} <\infty , \ \text{for all } \ p\in [1,\infty]. $
\end{proof}

For every $a\in (0,4)$ and $0<\varepsilon <3/8$, let us define   $M_\varepsilon := [4\varepsilon(1-\varepsilon)^2,1-\varepsilon]$  and the Markov chain ${Y^{a,b,\varepsilon}_{n+1}:=Y^\varepsilon_{n+1} = \omega_n Y_n^\varepsilon (1 - Y^\varepsilon_n)}$
absorbed at $\partial^\varepsilon = \R \setminus M_\varepsilon$, where $\{\omega_n\}_{n\in\mathbb N}$ is an i.i.d. sequence of random variables and $\omega_n\sim\mathrm{Unif}([a,b])$. Moreover, for every $\varepsilon \in (0,3/8),$  we denote the transition kernels and transfer operator for the absorbing Markov chain $Y_n^\varepsilon$ as, respectively, 
\begin{align}
\mathcal P_\varepsilon f(x) := \mathbbm 1_{M_\varepsilon}(x) \mathcal P (\mathbbm 1_{M_\varepsilon} f)(x) \text{ and }\mathcal L_\varepsilon f (x) := \mathbbm 1_{M_\varepsilon}(x) \mathcal L(\mathbbm 1_{M_\varepsilon} f)(x).\label{Peps}
\end{align}

In the next proposition we show  the existence of a sequence of positive real numbers $\{\varepsilon_i\}_{i\in\mathbb N}$ converging to $0$, such that for every $i\in\mathbb N$, the absorbing Markov chain $Y_n^{\varepsilon_i}$ admits a unique quasi-stationary measure $\mu_{\varepsilon_i}$ supported on $M_{\varepsilon_i}.$ Moreover, these measures will play an important role in constructing a non-trivial quasi-stationary measure for $Y_n$ on $M$.

\begin{proposicao}\label{TBK}
Let $(a,b) \in [1,4) \times (4,\infty)$  and $Y_n^{a,b,\varepsilon}$ be the Markov chain absorbed at $\partial^{\varepsilon}$ defined above. Then, there exists a sequence of positive numbers $\{\varepsilon_i\}_{i\in \mathbb N}$ converging to $0$ such that, for every $i\in\mathbb N$, the following items hold:
\begin{enumerate}
    \item[(a)]  $Y_n^{a,b,\varepsilon_i}$ admits a unique quasi-stationary measure $\mu_{a,b,\varepsilon_i} := \mu_\varepsilon$ on $M_\varepsilon$ with survival rate $\lambda_{\varepsilon_i}>0$;
    \item[(b)] there exists a continuous function $g_\varepsilon^{a,b}:= g_{\varepsilon_i} \in \mathcal C^0(M_{\varepsilon_i})$ such that $\mu_{\varepsilon_i}(\d x) = g_{\varepsilon_i}(x)\d x$; and
    \item[(c)] $\supp(\mu_{\varepsilon_i}) = M_{\varepsilon_i}.$
\end{enumerate}
\end{proposicao}

\begin{proof}

From \cite[Theorem B and Remark XIII/5]{Jakobson} there exists a sequence $\{r_i\}_{i\in\mathbb N} \subset [a,4)$ converging to $4$, such that for every $i\in\mathbb N$ the logistic map $f_{r_i} :[0,1]\to [0,1]$, $f_{r_i}(x)= r_i x(1-x)$ admits an invariant ergodic measure $\rho_{r_i} \ll \mathrm{Leb}$ and $\supp(\rho_{r_i}) = [f_{r_i}^2(1/2), f_{r_i}(1/2)]$.

Consider the sequence $\{\varepsilon_i = (4 - r_i)/4\}_{i\in\mathbb N}$. Combining \eqref{Peps} and Proposition \ref{uk1} we obtain that
\begin{align}
    \frac{\mathcal P_{\varepsilon_{i}}(x,\d y)}{\mathrm{Leb}(\d y)} = \frac{\mathbbm{1}_{M_{\varepsilon_i}}(y) } { { (b-a) x(1-a)}} \mathbbm{1}_{[a x(1-x), bx(1-x)]}(y),\ \text{for every }i\in\mathbb N. \label{RNder}
\end{align}

In the following, we show that for every $i\in\mathbb N,$ given $x\in M_{\varepsilon_i}$ and open interval $I \subset M_{\varepsilon_i} =[f_{r_i}^2(1/2), f_{r_i}(1/2)] ,$ there exists $n_0 = n_0(x,I) \in\mathbb N,$ such that $\mathcal P_{\varepsilon_i}^{n_0}(x, I) > 0.$ 

Consider the set $J := \{y\in M_{\varepsilon_i}; \omega x(1-x) =y\ \text{for some }\omega \in [a,b]\}.$ Since $J$ has non-empty, interior we obtain that $\rho_{r_i}(J) >0.$ Since $\rho_{r_i}$ is an invariant ergodic measure, there exists $\omega_0 \in [a,b]$ such that $y:=\omega_0 x(1-x) \in J$ and $n_1>0 \in\mathbb N$ such that $f^{n_1}_{r_i}(y) \in I$. 

Consider the natural number $n_0 \in \mathbb N$ and the continuous function ${F^{x,n_0}: [a,b]^{n_0} \to \mathbb R,}$ $F^{x,n_0}(c_1,\ldots, c_{n_0}) := f_{c_1} \circ f_{c_2} \circ \ldots f_{c_{n_0}}(x)$. From the last paragraph we obtain that $F^{x,n_0}(\omega_0,r_{\varepsilon_i},\ldots,r_{\varepsilon_i}) \in I$. Finally, since $F^{x,n_0}$ is a continuous function we obtain that
\begin{align}
    \mathcal P^{n_0}(x, I) &= \mathbb P [Y_{n_0}^{a,b,\varepsilon_i} \in I\mid X_0 = x]= \frac{1}{(b-a)^{n_0}}\mathrm{Leb}^{\otimes n_0}\left(\{p\in [a,b]^{n_0};\  F^{x,n_0}(p) \in I\}\right) >0. \label{Transitiv}
\end{align}

From \eqref{Peps} and \eqref{Transitiv} we conclude that \cite[Hypothesis (H)]{Me} is fulfilled and therefore items $(a),$ $(b)$ and $(c)$ follows directly from \cite[Theorem A]{Me}.

\end{proof}

% From \cite[Theorem A]{Me}, we obtain that for every $\varepsilon\in (0,1/2),$ there exists a unique quasi-stationary measure $\mu_{\varepsilon}^{a,b}:=\mu_{\varepsilon}$ on $[\varepsilon,1-\varepsilon]$ for $Y_n^\varepsilon$ with survival rate $\lambda_\varepsilon^{a,b} :=\lambda_\varepsilon$.  Moreover, for every $\varepsilon \in (0,1/2),$  we denote the transition kernels and transfer operator for the absorbing Markov chain $Y_n^\varepsilon$ as, respectively, 
% \begin{align*}
% \mathcal P_\varepsilon f(x) := \mathbbm 1_{M_\varepsilon}(x) \mathcal P (\mathbbm 1_{M_\varepsilon} f)(x) \text{ and }\mathcal L_\varepsilon f (x) := \mathbbm 1_{M_\varepsilon}(x) \mathcal L(\mathbbm 1_{M_\varepsilon} f)(x).
% \end{align*}

Observe that the family of measures given by the previous proposition $\{\mu_{\varepsilon_i}\}_{i\in \mathbb N}$ can be naturally extended on $[0,1]$ by imposing that $\mu_{\varepsilon_i}([0,1]\setminus M_{\varepsilon_i}) =0$ for every $i\in\mathbb N.$ In order to construct a quasi-stationary measure for the Markov process $Y_n$ on $[0,1]$ we use that $\{\mu_{\varepsilon_i}\}_{i \in \mathbb N}$ is precompact in the weak$^*$ of $\mathcal M([0,1]),$ i.e.
\begin{align}
    \bigcap_{i\in \mathbb N} \overline{ \{\mu_{\varepsilon_{k+i}}\}_{k\in\mathbb N}}^{\mathrm{w}^*\text{-}\mathcal M(M)} \neq \emptyset, \label{hv}
\end{align}
where $\mathrm{w}^*\text{-}\mathcal M(M),$ denotes the weak${}^*$ topology of $\mathcal M(M).$ 
% and 
% $$\mathcal L_\varepsilon f (x) = \mathbbm 1_{M_\varepsilon}(x) \mathcal L(\mathbbm 1_{M_\varepsilon} f)(x). $$

The proposition below shows that the elements of \eqref{hv} are natural candidates for quasi-stationary measures for $Y_n$ on $[0,1]$.

\begin{proposicao}
Assume that there exists a probability measure $\mu_{a,b} :=\mu$ on $M$, $\lambda>0,$  and subsequences $$\{\mu_{\delta_n}\}_{n\in\mathbb N} \subset \{\mu_{\varepsilon_n}\}_{n \in \mathbb N}, \ \{\lambda_{\delta_n}\}_{n\in\mathbb N} \subset \{\lambda_{\varepsilon_n}\}_{n\in \mathbb N}, $$
such that
$$\mu_{\delta_n} \to \mu, \ \text{in the weak}^{*}\text{-topology}\ \text{as }n\to\infty, $$
$$\lim_{n\to \infty} \lambda_{\delta_n} = \lambda \ \text{and }\lim_{n\to \infty }\delta_n = 0. $$

Then $\mu$ is a quasi-stationary measure for $Y_n$ on $[0,1].$
\label{xavitorras}
\end{proposicao}
\begin{proof}
Let 
$$E = \{x\in [0,1], \mu(\{x\}) >0\},$$
note that $E$ is, at most, countable. Consider the set
$$\mathcal A =\{I \in \mathcal B(M); \ I\text{ is an interval,}\ \overline{I}\subset (0,1)\ \text{and }\sup I,\inf I\not\in E\}. $$
It is clear that $\sigma(\mathcal A) = \mathcal B(M).$  Note that for every $I \in \mathcal A,$ there exists $n_0 = n_0(I),$ such that $$A \subset M_{\delta_n},\ \text{for all }\ n>n_0.$$
This implies that for every $n>n_0,$
\begin{align*}
 \int_M \mathcal P(x,I) \mu_{\delta_n}(\d x)  &= \int_{M_\varepsilon} \mathcal P(x,I) \mu_{\varepsilon_n}(\d x) = \lambda_{\delta_n} \mu_{\delta_n} (I).
\end{align*}
Since $\mathcal P(x,I)$ is  a continuous function, we obtain
$$ \int_M \mathcal P(x,I) \mu(\d x) =\lambda \mu(I) \ \text{for every }I\in\mathcal A.$$
Since $\frac{4-b}{b-a} \leq \mathcal P(x,M), \ \text{for all } \ x\in[0,1],$ it follows that $\lambda >0.$
 
Applying the monotone class theorem, we obtain that $\mu$ is a quasi-stationary measure for $Y_n$ on $[0,1]$.

\end{proof}

In light of Proposition \ref{xavitorras}, to construct a non-trivial quasi-stationary measure for $Y_n$ on $[0,1]$, it remains to show that 
\begin{align}
    \bigcap_{i\in \mathbb N} \overline{ \{\mu_{\varepsilon_{i+k}}\}_{i\in \mathbb N}}^{\mathrm{w}^*\text{-}\mathcal M(M)} \setminus \{\delta_0\} \neq \emptyset. \label{vh}
\end{align}

Note that for every $i\in\mathbb N$,  $\mu_{\varepsilon_i}(\d x) \ll \mathrm{Leb}(\d x)$. To show that \eqref{vh}  holds, we study the behaviour of the distributions of $\mu_{\varepsilon_i}$ with respect to the Lebesgue measure.

The definition below provides conditions on $a$ and $b$, which implies that \eqref{vh} holds (see Theorem \ref{fnx}).

\begin{definicao}\label{AdmissiblePair}
A pair $(a,b)\in \mathbb (0,4)\times (4,\infty)$ is called an \emph{admissible pair} if either

\begin{itemize}
    \item $a\geq 2;$ or
    \item for every $x\in [(4 a^2 - a^3)/16,a/4]$
\begin{align}
         0\leq \frac{1}{2}-\frac{1}{2} \sqrt{1-\frac{2}{b} \left(1- \sqrt{1-\frac{4
   x}{a}}\right)} \leq \frac{a}{4} \label{ineq1}
\end{align}
 and
      \begin{align}
       \frac{2 \left(\tanh ^{-1}\left(\sqrt{\frac{2 \sqrt{1-\frac{4 x}{b}}+b-2}{b}}\right)-\tanh
   ^{-1}\left(\sqrt{\frac{a+2 \sqrt{1-\frac{4 x}{b}}-2}{a}}\right)\right)}{2 \tanh ^{-1}\left(\sqrt{\frac{2 \sqrt{1-\frac{4 x}{a}}+b-2}{b}}\right)+\log
   \left(\frac{a}{4-a}\right)}\leq \frac{\sqrt{1-\frac{4
   x}{b}}}{\sqrt{1-\frac{4 x}{a}}}.\label{ineq2}
      \end{align}  
\end{itemize}
\end{definicao}
% \begin{remark}
% We mention that fixed $(a,b) \in (0,1)\times (4,\infty)$ it is relatively check if $(a,b)$ is an admissible pair or not. However, it is 
% \end{remark}

In Theorem \ref{T616} we show that if $(a,b)\in [1,4)\times (4,\infty)$ then $(a,b)$ is an admissible pair. Assuming that $(a,b)$ is an admissible pair it is possible to show that $Y_n^{a,b}$ admits a non-trivial quasi-stationary measure on $[0,1].$ To accomplish this goal, we need the following three technical lemmas.

\begin{lema}\label{milagre0}
Let $(a,b)$ be an admissible pair, with $a <2$ and $f:[0,1]\to \mathbb R$  be a function continuous by parts with a finite number of discontinuities, such that
\begin{enumerate}
    \item $0\leq f(x)$ for every $x\in [0,1];$
    \item $f$ is non-decreasing in the interval $[0,a/4];$ and
    \item $f$ is non-increasing in the interval $[a/4,1].$
\end{enumerate}
% and $ 2-\sqrt{(-2 + b)/b}<a.$

Then $\mathcal L f$ is a continuous function such that
\begin{enumerate}
    \item $0\leq \mathcal L f(x)$ for every $x\in [0,1];$
    \item $\mathcal L f$ is non-decreasing in the interval $[0,(4a^2 - a^3)/16];$ and
    \item $\mathcal L f$ is non-increasing in the interval $[a/4,1].$
\end{enumerate}
\label{d2}
\end{lema}

\begin{proof}
Recall that
\begin{align}
   \mathcal L f(x) =  \int_{\alpha_-(x)}^{\alpha_+(x)} \frac{f(y)}{(b-a)y(1-y)} \, \d y  - \int_{\beta_-(x\wedge a/4)}^{\beta_+(x\wedge a/4)} \frac{  f(y)}{(b-a)y(1-y)} \, \d y ,  \label{induc}
\end{align}
where
\begin{align*}
    \alpha_{\pm}(x) = \frac{1}{2} \pm \frac{1}{2}\sqrt{1 - \frac{4}{b}x}\ \text{and } \beta_{\pm}(x) = \frac{1}{2} \pm \frac{1}{2}\sqrt{1 - \frac{4}{a}x}.
\end{align*}

It is clear that $\mathcal L f$ is continuous and a non-negative function. Observe that $\mathcal L f$ is differentiable except for finitely many  points. In fact, the derivative of $\mathcal L f$ on the points that the derivative exists is given by
\begin{align}\label{eqn:Jetlag}
\frac{\d \mathcal L f}{\d x}(x) = -  \frac{ f(\alpha_+ (x)) + f(\alpha_-(x))}{ (b-a) x\sqrt{1-\frac{4}{b}x} } + \mathbbm 1_{\left[0,\frac{a}{4}\right]}(x)  \frac{f(\beta_+ (x)) + f(\beta_-(x))}{(b-a) x\sqrt{1- \frac{4}{a}x} }. 
\end{align}

Since for every $x\in [a/4,1],$
$$ \frac{\d \mathcal L f}{\d x}(x) = - \frac{1}{b-a}\frac{ f(\alpha_+ (x)) + f(\alpha_-(x))}{ x\sqrt{1-\frac{4}{b}x} } \leq 0,$$
if follows that $\mathcal L f$ is non-increasing in $[a/4,1].$ 

Observe that for every $x\in [0,a/4]$ we obtain
$$\frac{1}{(b-a) x\sqrt{1-\frac{4}{b}x} } <  \frac{1}{(b-a) x\sqrt{1- \frac{4}{a}x}},\ \text{and}\ \frac{a}{4} \leq \beta_+(x)\leq \alpha_+(x).$$
Since $f$ is non-increasing in $[a/4,1],$ we conclude that
$$-\frac{ f(\alpha_+(x))}{(b-a) x\sqrt{1-\frac{4}{b}x} }  + \frac{f(\beta_+(x))}{(b-a) x\sqrt{1- \frac{4}{a}x}}\geq 0. $$

% \begin{figure}[h!]
%     \centering
%     \includegraphics[scale=0.4]{graph4.pdf}
% \end{figure}

To finish the proof, it is enough to show that $f(\beta_-(x)) \geq f(\alpha_-(x))$ for every $x\in [0,(4a^2 - a^3)/16].$ Observe that since $f$ is non-decreasing on $[0,a/4]$, we obtain that for every $x\in [0,(4a^2 - a^3)/16]$ 
$$\beta_-(x) \leq \alpha_-(x) \leq a/4, $$
implying that
$$ f(\beta_-(x)) -f(\alpha_-(x)) \geq 0.$$
\end{proof}

% The following three lemmas study the properties of the distributions $\{g_\varepsilon\}.$

\begin{lema}
\label{milagre3} Let $(a,b)$ be an admissible pair and $\varepsilon \in (0,3/8)$ such that $[(4a^2 - a^3)/16, a/4]\subset {M_\varepsilon}$. Consider that sequence of functions  $\left\{\mathcal L_\varepsilon^n \mathbbm 1_{M_\varepsilon} \right\}_{n\in\mathbb N}$, then for every $n\in\mathbb N$ the following assertions hold:

\begin{enumerate}
    \item $0\leq \mathcal L^n_\varepsilon \mathbbm 1_{M_\varepsilon}(x)$ for every $x\in [0,1];$
    \item $ \mathcal L_\varepsilon ^n \mathbbm 1_{M_\varepsilon}(x)$ in non-decreasing in the interval $[0,(4a^2-a^3)/16];$ and
    \item $ \mathcal L_\varepsilon ^n \mathbbm 1_{M_\varepsilon}(x)$ is non-increasing in the interval $[a/4,1].$
\end{enumerate}
\end{lema}

\begin{proof}

Recall that for every $\varepsilon \in (0,3/8)$ and $f\in\mathcal C^0(M_{\varepsilon}),$ $\mathcal L_{\varepsilon} f   = \mathbbm 1_{M_{\varepsilon}} \mathcal L (\mathbbm 1_{M_{\varepsilon}} f ).$  We divide the proof into two steps.

\begin{step}
\it We show the result for the case that $a\geq 2.$
\end{step}
We show the above result by induction on $n$. The case $n=0$ is immediately verified. Suppose that items $(1),$ $(2)$ and $(3)$ are true for $\mathcal L_{\varepsilon}^{n} \mathbbm 1_{M_{\varepsilon}},$ we will show that the same is true for $\mathcal L^{n+1}_{\varepsilon} \mathbbm 1_{M_{\varepsilon}}.$

Item $(1)$ is trivially fulfilled since $\mathcal L_\varepsilon$ is a positive operator. Additionally items $(2)$ and $(3)$ follow from equation \eqref{eqn:Jetlag}  and realising that for every $(a,b)\in [2,4)\times (4,\infty)$
$$\alpha_-(x) \leq \beta_-(x) \leq \frac{4a^2 - a^3}{16} \leq \frac{a}{4}\leq \beta_+(x) \leq \alpha_+(x),\ \text{for every }x\in\left[0,\frac{4 a^2 - a^3}{16}\right].$$

This proves Step 1.

\begin{step}
We show that if $(a,b)$ is an admissible pair and $a\in(0,2)$ then:
\begin{enumerate}
    \item $0\leq \mathcal L^n_\varepsilon \mathbbm 1_{M_\varepsilon}(x)$ for every $x\in [0,1];$
    \item $ \mathcal L_\varepsilon^n \mathbbm 1_{M_\varepsilon}(x)$ in non-decreasing in the interval $[0,a/4];$ and
    \item $ \mathcal L^n_\varepsilon \mathbbm 1_{M_\varepsilon}(x)$ is non-increasing in the interval $[a/4,1]$.
\end{enumerate}
\end{step}
We will prove that the above items hold by strong induction on $n$. The cases $n=0$ and $n = 1$, the computations can explicitly be done and such a conclusion is achieved.

Now, suppose that the conclusions of Step $2$ are true for $$\mathbbm 1_{M_\varepsilon}, \mathcal L_\varepsilon^1 \mathbbm 1_{M_\varepsilon}, \ldots, \mathcal L_\varepsilon^n \mathbbm 1_{M_\varepsilon}\ \text{with }n\geq 1$$ with and we will show that it is also true for $\mathcal L_\varepsilon^{n+1} \mathbbm 1_{M_\varepsilon}.$ 

From Lemma \ref{d2} it follows that
\begin{enumerate}
    \item $0\leq \mathcal L^{n+1}_\varepsilon \mathbbm 1_{M_\varepsilon}(x)$ for every $x\in [0,1];$
    \item $\mathcal L^{n+1}_{\varepsilon} \mathbbm 1_{M_\varepsilon}$ in non-decreasing in the interval $[0,(4a^2 - a^3)/16];$ and
    \item $\mathcal L^{n+1}_\varepsilon \mathbbm 1_{M_{\varepsilon}}$ is non-increasing in the interval $[a/4,1].$
\end{enumerate}

It remains to show that $\mathcal L^{n+1}_\varepsilon {M_{\varepsilon}} $ is non-decreasing in $[(4a^2 - a^3)/16,a/4].$  From the proof of the previous theorem, it is enough to show that
\begin{align}
\frac{\mathcal L^n_\varepsilon \mathbbm 1_{M_{\varepsilon}} (\alpha_- (x))}{  \sqrt{1- \frac{4x}{b} }}\leq \frac{\mathcal L_\varepsilon^{n}\mathbbm 1_{M_{\varepsilon}} (\beta_-(x))}{ \sqrt{1 - \frac{4x}{a}}}\ \text{for every } x\in \left[\frac{4a^2 - a^3}{16},\frac{a}{4}\right].\label{verdade}
\end{align}

Observe that  $$\alpha_-(x) <\frac{a}{4} < \beta_-(x) \ \text{for every } x\in \left[\frac{4a^2 - a^3}{16},\frac{a}{4}\right]. $$
% \begin{figure}[h!]
%     \centering
%     \includegraphics[scale=0.3]{graph1.pdf}
% \end{figure}

Therefore,
\begin{align}
   \mathcal L_\varepsilon^{n}\mathbbm 1_{M_{\varepsilon}} (\beta_-(x)) = \int_{\alpha_-\circ \beta_-(x)}^{\alpha_+ \circ \beta_-(x)} \frac{ \mathcal L_\varepsilon^{n-1}\mathbbm 1_{M_{\varepsilon}}(y)}{(b-a) y (1-y)} \, \d y \label{vincent1}
\end{align} 
and
\begin{align*}
\mathcal L_\varepsilon^{n}\mathbbm 1_{M_{\varepsilon}} (\alpha_-(x)) &= \int_{\alpha_-\circ \alpha_-(x)}^{\alpha_+ \circ \alpha_-(x)} \frac{ \mathcal L_\varepsilon^{n-1}\mathbbm 1_{M_{\varepsilon}}(y)}{(b-a) y (1-y)} \, \d y - \int_{\beta_-\circ \alpha_-(x)}^{\beta_+ \circ \alpha_-(x)} \frac{ \mathcal L_\varepsilon^{n-1}\mathbbm 1_{M_{\varepsilon}}(y)}{(b-a) y (1-y)} \, \d y \\
&=  \int_{\alpha_-\circ \alpha_-(x)}^{\beta_- \circ \alpha_-(x)} \frac{ \mathcal L_\varepsilon^{n-1}\mathbbm 1_{M_\varepsilon}(y)}{(b-a) y (1-y)} \, \d y + \int_{\beta_+\circ \alpha_-(x)}^{\alpha_+ \circ \alpha_-(x)} \frac{ \mathcal L_\varepsilon^{n-1}\mathbbm 1_{M_\varepsilon}(y)}{(b-a) y (1-y)} \, \d y. 
\end{align*}

Since $(a,b)$ is an admissible pair, \eqref{ineq1} implies that for every  $x\in \left[(4a^2 - a^3)/16,a/4\right]$
$$  \beta_-\circ\alpha_-(x) < \alpha_- \circ \beta_- (x)  <\frac{a}{4}<  \alpha_+ \circ \beta_-(x) \leq \beta_+ \circ \alpha_-(x) . $$

This implies that
\begin{align*}
    \mathcal L^n_\varepsilon\mathbbm 1_{M_{\varepsilon}}(\alpha_-(x)) \leq &\frac{ \mathcal L_\varepsilon^{n-1}\mathbbm 1_{M_{\varepsilon}}(\beta_-\circ \alpha_-(x))}{b-a}\int_{\alpha_-\circ \alpha_-(x)}^{\beta_- \circ \alpha_-(x)} \frac{ 1}{ y (1-y)} \, \d y \\ 
    &+\frac{\mathcal L_\varepsilon^{n-1}\mathbbm 1_{M_{\varepsilon}}(\beta_+\circ \alpha_- (x) )}{b-a}\int_{\beta_+\circ \alpha_-(x)}^{\alpha_+ \circ \alpha_-(x)} \frac{ 1}{ y (1-y)}\, \d y \nonumber\\
    \leq&  \frac{  \mathcal L_\varepsilon^{n-1}\mathbbm 1_{M_{\varepsilon}}(\beta_-\circ \alpha_-(x))  + \mathcal L_\varepsilon^{n-1}\mathbbm 1_{M_{\varepsilon}}\left( \beta_+\circ \alpha_- (x)\right)}{b-a} I_1^{(a,b)}(x)
\end{align*}
where
$$I_1^{(a,b)}(x):=2  \left(\tanh ^{-1}\left(\sqrt{\frac{2 \sqrt{1-\frac{4 x}{b}}+b-2}{b}}\right)-\tanh
   ^{-1}\left(\sqrt{\frac{a+2 \sqrt{1-\frac{4 x}{b}}-2}{a}}\right)\right).$$
On the other hand, from the induction hypothesis and \eqref{vincent1}
\begin{align*}
    \mathcal L^{n}_\varepsilon\mathbbm 1_{M_{\varepsilon}} (\beta_-(x)) \geq& \frac{\mathcal L_\varepsilon^{n-1}\mathbbm 1_{M_{\varepsilon}}(\beta_-\circ \alpha_-(x))}{b-a}\int_{\alpha_-\circ \beta_-(x)}^{a/4} \frac{ 1}{y (1-y)} \, \d y  \\
    +&\frac{\mathcal L_\varepsilon^{n-1}\mathbbm 1_{M_\varepsilon}(\beta_+\circ \alpha_- (x) )}{b-a}\int_{a/4}^{\alpha_+ \circ \beta_-(x)} \frac{ 1}{ y (1-y)}\, \d y\nonumber.\\
    \geq&  \frac{\mathcal L_\varepsilon^{n-1}\mathbbm 1_{M_\varepsilon}(\beta_-\circ \alpha_-(x))+ \mathcal L^{n-1}_\varepsilon\mathbbm 1_{M_\varepsilon}\left( \beta_+\circ \alpha_- (x)\right)}{b-a}I_2^{(a,b)}(x), 
\end{align*}
where
$$I_2^{(a,b)}(x) :=\left( 2\tanh ^{-1}\left(\sqrt{\frac{2 \sqrt{1-\frac{4 x}{a}}+b-2}{b}}\right)+\log
   \left(\frac{a}{4-a}\right)\right).$$

Combining the above three equations, \eqref{ineq2} and using the definition of admissible pair we obtain that \eqref{verdade} holds. This proves Step $2$. 

Observe the above two steps implies the proof of the lemma.

\end{proof}

Recall from Proposition \ref{TBK}, for every $i\in\mathbb N$
$$g_{\varepsilon_i}^{a,b}:= g_{\varepsilon_i} = \frac{\mu_{\varepsilon_i}(\d x) }{\mathrm{Leb}(\d x)}  \in L^1 ([0,1],\mathrm{Leb}),$$
where we set $g_{\varepsilon_i}(x) =0$ for every $x\in M\setminus M_{\varepsilon_i}.$

\begin{lema}\label{structure}
Let $(a,b)$ be an admissible pair. Then, for every $i\in\mathbb N$
\begin{enumerate}
    \item $0\leq g_{\varepsilon_i} (x)$ for every $x\in [0,1];$
    \item $ g_{\varepsilon_i}(x)$ in non-decreasing in the interval $[0,(4a^2 - a^3)/16];$ and
    \item $ g_{\varepsilon_i} (x)$ is non-increasing in the interval $[a/4,1].$
\end{enumerate}

\label{ge}

\end{lema}

\begin{proof}

Recall that for every $i\in \mathbb N$ and $f\in\mathcal C^0(M_{\varepsilon_i}),$ $\mathcal L_{\varepsilon_i}f   = \mathbbm 1_{M_{\varepsilon_i}} \mathcal L (\mathbbm 1_{M_{\varepsilon_i}} f ).$ Observe that if $(a,b)$ is an admissible pair and $i\in\mathbb N$ then  $\mathcal L_{\varepsilon_i} :\mathcal C^0(M_{\varepsilon_i} )\to \mathcal C^0(M_{\varepsilon_i} )$ is an irreducible compact operator. Moreover, it is readily verified that $\mathcal L_{\varepsilon_i} $ admits a single eigenvalue in its peripheral spectrum, implying that
$$\left\|\frac{1}{\lambda_{\varepsilon_i}^n}\mathcal L^n_{{\varepsilon_i} }\mathbbm 1_M  -  \alpha_{\varepsilon_i}  g_{\varepsilon_i}  \right\|_{\mathcal C^0(M_{\varepsilon_i} )} \xrightarrow{n\to\infty} 0, $$
for some $\alpha_{\varepsilon_i}  >0.$

The Lemma follows directly from the above equation in combination with Lemma \ref{milagre3}.

\end{proof}

Combining  Lemmas \ref{milagre0}, \ref{milagre3}, and \ref{structure} we obtain the following result.

\begin{teorema}
Let $(a,b)$ be an admissible pair. Then the absorbing Markov chain $Y_n^{a,b}$ admits a quasi-stationary measure $\mu$ on $[0,1]$ different from $\delta_0.$ \label{fnx}
\end{teorema}
\begin{proof}
For every $i\in \mathbb N$, let  $\mu_{\varepsilon_i}(\d x) = g_{\varepsilon_i} (x) \d x$ to be the unique quasi-stationary measures for $Y_n^{\varepsilon_i}$ on $M_{\varepsilon_i}$ given by Proposition \ref{TBK} and extend it to $[0,1]$ in a way that $\mu_{\varepsilon_i}( M\setminus M_{\varepsilon_i} ) = 0. $

Since $\mathcal M_1 ([0,1])$ is sequentially compact in the weak$^*$ topology, we can assume without loss of generality  (passing to a subsequence if necessary) that the sequence of real numbers $\{\varepsilon_i\}_{i\in\mathbb N}$ is such that
$\lim_{i\to \infty} \varepsilon_i = 0,  \mu_{\varepsilon_i} \to \mu \ \text{in the weak}^{*}\ \text{topology}$
and $\lim_{i\to \infty} \lambda_{\varepsilon_i} = \lambda \geq \frac{4-a}{b-a}.$

From Proposition \ref{xavitorras}, the probability measure $\mu$ is a quasi-stationary measure of $\mathcal P$. It remains to show that $\mu \neq \delta_0.$ Suppose by contradiction that $\mu = \delta_0.$ Then $\lim_{i\to\infty} \mu_{\varepsilon_i}([0,(4a^2 -a^3)/32]) = 1. $ However, from Lemma \ref{ge}, it follows that $${\mu_{\varepsilon_i}([0,(4a^2 - a^3)/32]) \leq \mu_{\varepsilon_i}([(4a^2 -a^3)/32,4a^2 -a^3)/16])}\text{ for every }i\in\mathbb N.$$ Taking the limit as $i\to\infty$ we obtain that $1 \leq \mu([4a^2 -a^3)/32,(4a^2 -a^3)/16])$, which is a contradiction. Implying that $\mu \neq \delta_0.$ 

\end{proof}

\begin{remark}\label{Problem}
Observe that without assuming that $(a,b)$ is an admissible pair, the inductive step presented in Step 2 of Lemma \ref{milagre3} no longer holds. Without this lemma, the core argument in the proof Theorem \ref{fnx} cannot be applied, and the existence of a non-trivial quasi-stationary measure for $Y_n^{a,b}$ becomes unclear.
\end{remark}

From now on we define $\mu_{a,b} = \mu$ as a non-trivial quasi-stationary measure for $Y_n$ on $[0,1]$ and $\lambda_{a,b} =\lambda$ its associated survival rate (given by Theorem \ref{fnx}). The next proposition shows that $\mu$ is absolutely continuous to the Lebesgue measure.

\begin{proposicao}
Let $(a,b)$ be an admissible pair. Then  $\mu \ll \mathrm{Leb}(\d x),$ and $0<\lambda < 1.$
\end{proposicao}
\begin{proof}
We can decompose $\mu (\d x) = \mu(\{0\})\delta_0(\d x) + \mu'(\d x) + \mu(\{1\})\delta_1(\d x).$

Since $\delta_0 \neq \mu$, we obtain that $\mu(\{0\}) \neq 1.$ Observe that
\begin{align*}
    \lambda \mu(\{0\}) \delta_0(\d x) + \lambda \mu'(\d x) + \lambda \mu(\{1\}) \delta_1(\d x) &= \lambda \mu(\d x)\\
    &= \mathcal P^*(\mu)(\d x)\\
    &= (\mu(\{1\}) + \mu\{0\}) \delta_0(\d x)  + \mathcal P^* (\mu').
\end{align*}

Since $\mathcal P^*(\mu') \ll \mathrm{Leb}(\d x)$, it follows that $\mathcal P^*(\mu')(\{1\}) =0.$ Implying that
$$\mu(\{1\}) = 0,$$
and
$$\lambda \mu(\{0\}) = \mu(\{0\}). $$

We claim that $\lambda < 1.$ Suppose by contradiction the opposite that $\lambda =1$, then 
$$\mu = {(\mathcal P^*)}^n \mu =  \mu(\{0\})\delta_1(\d x) + ({\mathcal P^*})^n \mu'.   $$
Since  $\mathcal P^n1 (x) \to \mathbbm 1_{\{0\}\sqcup \{1\}}(x) $
 pointwise as $n\to \infty.$ We have, by the Lebesgue-dominated convergence theorem
$$\mu = \lim_{n\to\infty}{(\mathcal P^*)}^n \mu =  \mu\{0\}\delta_1(\d x) + \lim_{n\to\infty}{\mathcal P^*}^n \mu' =   \mu(\{0\}) \delta_0, $$
which is contradiction since $\mu(\{0\}) \neq 1.$

Implying that $\lambda < 1$ and therefore $\mu(\{0\}) = 0.$  Therefore $\mu (\{0\}\cup \{1\}) = 0$, and 
$$ \lambda \mu = \mathcal P^* \mu \ll \mathrm{Leb}(\d x). $$
\end{proof}
From now on, we define
$$\frac{\mu_{a,b} (\d x)}{\mathrm{Leb}(\d x)} =: g^{a,b}= g \in L^1([0,1],\mu). $$
The next result summarises the properties of $g.$

% \begin{proposicao}
% Let $\{\mu_{\varepsilon_n}\}_{n\in\mathbb N}$ be the sequence of measures such that  
% $$g_{\varepsilon_n}\, \d x = \mu_{\varepsilon_n}(\d x) \to \mu \ \text{in the weak}^* \ \text{topology}, $$
% then 
% $$g_{\varepsilon_n} \to g\ \text{in }L^{1}([0,1]). $$
% \end{proposicao}

% \begin{proof}
% note that for every interval $(a,b) \subset [0,1],$ we have that
% $$\lim_{n\to \infty} \int_a^b g_{\varepsilon_n}(y) \d y = \int_a^b g(y) \ \d y, $$
% by a monotone class theorem we obtainthat
% $$ \lim_{n\to \infty} \int_B g_{\varepsilon_n}(y) \d y = \int_B g(y) \ \d y \ \dieKK$$
% \end{proof}

\begin{proposicao}\label{imperial}
Let $(a,b)$ be an admissible pair. Then the function $g$ fulfils the following properties
\begin{enumerate}
    \item[(i)] $g \in \mathcal C^0(M);$
    \item[(ii)] $g$ is non-decreasing in the interval $[0,4a^2 -a^3)/16];$
    \item[(iii)] $g$ is non-increasing in the interval $[a/4,1]$;
    \item[(iv)] there exists $k>0$ such that $k<g(x)$ for every $x\in M.$ 
\end{enumerate}
\end{proposicao}
\begin{proof}
We divide this proof into $3$ steps.

\begin{step} \it We show that $g(x)>0,$ for every $x\in(0,1].$ 
\end{step}
Suppose that there exists $x \in (0,1]$ such that $g(x) = 0.$ Therefore
$$0 = \lambda g(x) = \int_{\alpha_-(x)}^{\alpha_+(x)}\frac{g(y)}{(b-a)y(1-y)}\, \d y  - \int_{\beta_-(x\wedge a/4)}^{\beta_+(x\wedge a/4)}\frac{g(y)}{(b-a)y(1-y)}\, \d y.   $$

This implies that 
$$ g(y) = 0,\ \text{for all } y\in I_1:= [\alpha_-(x),\beta_-(x\wedge a/4)]\cup [\beta_+(x\wedge a/4),\alpha_+(x)] \subset (0,1).$$

Let $x_0 \in \supp(\mu)\cap (0,1).$ By the same arguments presented in the proof of Proposition \ref{TBK}, we can show that there exist $n_0 = n_0(x,I_1)$ such that $\mathcal P^{n_0}(x_0, I_1)>0.$ Since $\mathcal P^{n_0}(x_0,I_1)$ is a continuous function, there exists an open neighbourhood $B \subset (0,1)$ of $x$, such that $$\inf_{y\in B} \mathcal P^{n_0}(y,I_1) \geq \frac{1}{2} \mathcal P^{n_0}(x_0,I_1) >0.$$ 
Therefore,
$$0 = \mu(I_1) = \frac{1}{\lambda^n}\int_M \mathcal P^{n_0}(y,I)g(y)\d x \geq \frac{\mathcal P^{n_0}(x,I_1)}{2}\mu(B)>0, $$
which is a contradiction. Therefore, $g(x) >0$ for every $x\in(0,1].$

\begin{step}\it
We show $(i),(ii)$ and $(iii).$ \label{structure1}
\end{step}

Recall that for every $i\in\mathbb N$, $ g_{\varepsilon_i}(x) = \frac{1}{ \lambda_{\varepsilon_i}}\mathbbm 1_{M_{\varepsilon_i}} \mathcal L g_{\varepsilon_i}(x).$ This observation, combined with Theorem \ref{structure} and Lemma \ref{d2} implies that
$$\|\cL g_{\varepsilon_i}\|_{L^\infty} = \sup_{y\in  \left[\frac{4a^2 - a^3}{16},\frac{a}{4}\right]}\cL g_{\varepsilon_i}(y), \ \text{for every }i\in\mathbb N. $$

Let 
$$ J := \bigcup_{x\in\left[\frac{4a^2 - a^3}{16},\frac{a}{4}\right]} \left( [\alpha_-(x),\beta_-(x) \land a/4 ]\cup [\beta_-(x) \land a/4, \alpha_+(x) \land a/4]\right)\subset (0,1),$$
and observe that $J$ is a compact set. Finally
\begin{align}
0&\leq g_{\varepsilon_i}(x) \leq  \frac{1}{\lambda_{\varepsilon_i}} \cL g_{\varepsilon_i}(x) \leq \frac{1}{\lambda_\varepsilon}\sup_{y\in \left[\frac{4a^2 - a^3}{16},\frac{a}{4}\right]}\cL g_{\varepsilon_i}(y)\nonumber\leq  \frac{1}{\lambda_{\varepsilon_i}}\int_{J}\frac{g_{\varepsilon_i} (y)}{(b-a)y(1-y)}\, \d y\\  &\leq \sup_{y\in J}\frac{1}{4 y (1-y)} \sup_{i\in\mathbb N} \frac{1}{\lambda_{\varepsilon_i}}=: C  <\infty \label{bound1}.
\end{align}

Therefore, we obtain an uniform bound for  $\{ g_{\varepsilon_i}\}_{i\in\mathbb N}$ on $L^{\infty}(M),$ for $n$ big enough. For every $\delta>0,$ consider the map
\begin{align*}
    T_\delta: L^{\infty}(M) &\to \mathcal C^0([\delta, 1-\delta])\\
    f&\to \mathbbm 1_{[\delta,1-\delta]} \mathcal L f.
\end{align*}
From the Arzelà–Ascoli theorem, it is readily verified that  $T_\delta$ is a compact operator for every $0 <\delta <1/2$. From \eqref{bound1}, we obtain that there exists a subsequence $\{\mathcal L g_{\varepsilon_{i_n}}\}_{n\in\mathbb N} \subset \{\mathcal L g_{\varepsilon_{i}}\}_{i\in\mathbb N} $ and $f_\delta \in \mathcal C^0([\delta,1-\delta])$ such that
$$\lim_{n\to \infty} \| T_\delta \mathcal L g_{\varepsilon_{i_n}} - f_\delta\|_{\infty} = 0.$$

Choosing an interval $I_\delta \subset [\delta,1-\delta]$, observe that
$$\mathcal P(x,I_\delta)\ \text{is continuous on }x\in [0,1], $$
it follows that
\begin{align*} 
\int_{I_\delta}  f_{\delta}(y) \d x = \lim_{n\to\infty } \int_{I_\delta} \mathcal L g_{\varepsilon_{i_n}}(x)\d x =  \lim_{n\to\infty } \int_{0}^{1} \mathcal P (x,I_\delta) g_{\varepsilon_{i_n}}(x)\d x =  \int_{I_\delta} \lambda g (x)\d x.
\end{align*}
Since $I_\delta$ is an arbitrary interval subset of $[\delta,1-\delta]$ we obtain $f_\delta =\lambda \left.g\right|_{[\delta,1-\delta]}$. Using that for every subsequence of $\{\mathbbm 1_{[\delta,1-\delta]} \mathcal L g_{\varepsilon_i}\}_{i\in\mathbb N}$ there exists subsubsequence converging to $ \mathbbm 1_{[\delta,1-\delta]} \lambda g$, we obtain that 
$$\lim_{i\to\infty}\|\mathbbm 1_{[\delta,1-\delta]} (\cL g_{\varepsilon_i} - \lambda g)\|_{\infty}  = 0.$$

Since $\{g_{\varepsilon_i}\}_{i\in\mathbb N}$ is bounded $L^\infty (M,\mu)$ and $g$ lying in $L^1(M,\mu),$ the above equation implies that
$$ \cL g_{\varepsilon_i} \to \lambda g\ \text{in }L^1([0,1]).$$

Thus, there exists a subsequence $\{g_{\varepsilon_{n_i}}\}_{i\in\mathbb N} \subset \{g_{\varepsilon_{n}}\}_{n\in\mathbb N}, $
such that
\begin{eqnarray}
\lim_{i\to\infty} \mathcal L g_{\varepsilon_{n_i}}  = \lambda g \ \mu\text{-almost surely}.\label{bs}
\end{eqnarray} 

Therefore for $\mu$-almost every $x\in M$
\begin{align*}
0\leq  g(x) &\leq \frac{1}{\lambda}\lim_{i\to\infty}\cL g_{\varepsilon_{n_i}}(x)\leq \frac{C}{\lambda },
\end{align*}
which implies that $g$ is $L^\infty([0,1])$. Since for every $i\in\mathbb N,$
\begin{enumerate}
    \item $ \cL g_{\varepsilon_{n_i}}$ in non-decreasing in the interval $[0,(4a^2 - a^3)/16];$ and
    \item $ \cL g_{\varepsilon_{n_i}}$  is non-increasing in the interval $[a/4,1].$\label{structure2}
\end{enumerate}
From \ref{bs} and the continuity of $g$ on $(0,1]$, we obtain that
\begin{enumerate}
    \item $g$ in non-decreasing in the interval $[0,(4a^2 - a^3)/16];$ and
    \item $g$  is non-increasing in the interval $[a/4,1].$\label{structure3}
\end{enumerate}
The proof is finished observing that $g\in \mathcal C^0([0,1])$ when imposing $g(0) := \inf_{x\in(0,a/4)} g(x). $
\begin{step} We show (iv). \label{g0}
\end{step}
Observe that in virtue of Step \ref{structure1}, it is enough to show that $g(0)>0$. 
Since $g$ is continuous, it follows that
$$\lim_{\varepsilon \to 0} \frac{1}{\varepsilon}\int_0^\varepsilon g(y) \d y = g(0). $$

Since $g(x)\d x$ is a quasi-stationary measure of $Y_n$ on $[0,1]$, we obtain that
\begin{align*}
    \int_0^\varepsilon g(y)\d y = \frac{1}{\lambda} \int_0^1 \mathcal P(y,[0,\varepsilon]) g(y) \d y.
\end{align*}
It is clear that 
$$\mathcal P(x,[0,\varepsilon]) = 1\ \text{for every } x \in [\alpha_+(\varepsilon), 1]\subset [a/4,1]. $$

Since $g$ is decreasing in $[a/4,1]$ and $g(1) >0$ it follows that
$$    \int_0^\varepsilon g(y)\d y = \frac{1}{\lambda} \int_0^1 \mathcal P(y,[0,\varepsilon]) g(y) \d y \geq \frac{g(1)}{\lambda} \left(1- \alpha_+(\varepsilon)\right) = \frac{g(1)}{\lambda} \left(\frac{1}{2}-\frac{1}{2} \sqrt{1-\frac{4 \varepsilon }{b}}\right).$$

Finally
\begin{align*}
    g(0) &= \lim_{\varepsilon \to 0} \frac{1}{\varepsilon} \int_0^\varepsilon g(y)\d y \geq  \lim_{\varepsilon \to 0} \frac{1}{\varepsilon} \frac{g(1)}{\lambda} \left(\frac{1}{2}-\frac{1}{2} \sqrt{1-\frac{4 \varepsilon }{b}}\right) = \frac{g(1)}{b\lambda} >0.
\end{align*}

Combining Steps 1-3, we conclude the proof of the theorem.

\end{proof}

To apply Theorem \ref{teorema4}, we need to show that $\mathcal P$ admits an eigenfunction lying $L^1([0,1],\mathrm{Leb}).$ To do this, consider the operator
\begin{align*}
    T: \mathcal C^0([0,1])&\to \mathcal C^0([0,1])\\
    f&\mapsto \frac{\mathcal L(g f)}{\lambda g},
\end{align*}
it is clear that $T$ is a Markov operator i.e.
\begin{enumerate}
    \item $T: \mathcal C^0([0,1])\to \mathcal C^0([0,1])$ is a bounded positive linear operator.
    \item $T1 =1.$
\end{enumerate}

\begin{proposicao}
Let $(a,b)$ be an admissible pair. Then there exists a probability $\nu \ll \mathrm{Leb}$ such that $\nu$ is a fixed point of the operator $T^*: \mathcal M(M) \to \mathcal M(M)$. \label{vidaboa}
\end{proposicao}
\begin{proof}
Since $T$ is a Markov operator, it is well known that there exists a probability measure $\nu$ such that $T^*\nu = \nu$ (see \cite[Chapter 10]{OPAET}).

Let us decompose $\nu$ as
$$\nu = \alpha_1 \delta_0 + \alpha_2 \nu' + \alpha_3 \delta_1,  $$
where $\nu\in \mathcal M_1(M)$ and $\nu'(\{0\} \cup \{1\}) = 0.$

Since
\begin{align*}
    \mathcal L (fg)(0) &= \frac{1}{b-a}\lim_{x\to 0} \left(\int_{\alpha_-(x)}^{\beta_-(x)} \frac{f(y)g(y)}{ y(1-y)}\, \d y + \int_{\beta_+(x)}^{\alpha_+(x)} \frac{f(y)g(y)}{ y(1-y)}\, \d y  \right) \\
    &= \frac{\log(b/a)}{b-a}f(0)g(0) + \frac{\log(b/a)}{b-a}f(1)g(1), 
\end{align*}
we obtain that
\begin{align}
    Tf(0) =\frac{\log(b/a)}{(b-a) \lambda }f(0) + \frac{\log(b/a)}{(b-a)\lambda} \frac{g(1)}{g(0)}f(1) \label{m3}.
\end{align}
From a similar computation, we obtain that
\begin{align}
Tf(1) = \frac{1}{ \lambda g(1)}\int_{\alpha_-(1)}^{\alpha_+(1)} f(x) g(x) \d x. \label{m22}
\end{align}

Note that given $A\in \mathscr B([0,1])$ such that  $\mathrm{Leb}(A)= 0$ and $A \subset [\delta, 1-\delta]$ for some $\delta>0,$ then $T(\mathbbm 1_A) = 0.$ This implies that
$$T^*\nu' (A) = \int_0^1 T^*\mathbbm 1_A(x) \nu'(\d x) =0,  $$
since $\nu'(\{0\} \cup \{1\}) = 0$, we obtain
\begin{align}
T^*\nu'(\d x) \ll \mathrm{Leb}(\d x). \label{m11}     
\end{align}

Combining $T^*\nu = \nu$, and equations \eqref{m3}, \eqref{m22} and \eqref{m11} we obtain that $\nu' \ll \mathrm{Leb}(\d x). $

Let $\{f_n\}_{n\in\mathbb N} \in \mathcal C^0(M)$ be a sequence of continuous functions such that
\begin{enumerate}
    \item $0\leq f_n (x)\leq 1$, for every $n\in\mathbb N$ and $x\in [0,1].$
    \item $f_n (1) =1;$ and
    \item $f_n(x) =0$ for every $x\in[0,1-1/n].$
\end{enumerate}

Since $T^*\nu = \nu$ and $f_n$ is continuous, it follows that
\begin{align}
    \int_M f_n(x) \nu(\d x)  = \int_M Tf_n(x) \nu(\d x) \ \text{for every }n\in\mathbb N. \label{bernat1}
\end{align}

The left-hand side of \eqref{bernat1} is equal to 
$$  \int_M f_n(x) \nu(\d x)  = \alpha_2 \int f_n(x) \nu(\d x) + \alpha_3  f_n(1) =   \alpha_2 \int f_n(x) \nu(\d x) + \alpha_3, $$
and the right-hand side of \eqref{bernat1} is equal to 
\begin{align*}
    \int_M f_n(x) T^*\nu(\d x)  =&  \alpha_1\left( \frac{\log(b/a)}{(b-a)\lambda }f_n(0) + \frac{\log(b/a) g(1)}{(b-a) g(0)\lambda}f_n(1)\right) + \alpha_2 \int_0^1 f_n(x)T^*\nu(\d x)  \\
    & + \alpha_3 \frac{1}{\lambda g(1)}\int_{\alpha_-(1)}^{\alpha_+(1)} f_n(x) g(x) \d x\\
    =& \alpha_1 \frac{\log(b/a)}{b-a} \frac{ g(1)}{\lambda g(0)} +  \alpha_2 \int_0^1 f_n(x)T^*\nu(\d x)+  \alpha_3 \frac{1}{\lambda g(1)}\int_{\alpha_-(1)}^{\alpha_+(1)} f_n(x) g(x) \d x.
\end{align*}  

Taking the limit as $n\to \infty$ in \eqref{bernat1} we obtain that
$$\alpha_3 = \alpha_1 \frac{\log(b/a)}{b-a} \frac{ g(1)}{\lambda g(0)}. $$

Repeating the same argument with the sequence $\{f_n(1-x)\}_{n\in\mathbb N} \subset \mathcal C^0([0,1])$ we obtain that
$$ \alpha_1 =  \frac{\log(b/a)}{(b-a)\lambda} \alpha_1.   $$

If $\alpha_1=0$ then $\alpha_3 =0$ and the proof is finished. Suppose by contradiction that $\alpha_1 >0$, the above equation shows that $ 1=  \log(b/a)\lambda^{-1}(b-a)^{-1}.$ On the other hand, we obtain that
$$g(0)= \frac{1}{\lambda}\mathcal L g(0) = \frac{1}{\lambda}\frac{\log(b/a)}{b-a}g(0) +\frac{1}{\lambda}\frac{\log(b/a)}{b-a}g(1) = g(0) +\frac{1}{\lambda}\frac{\log(b/a)}{b-a}g(1),  $$
therefore $g(1) =0,$ contradicting  Proposition \ref{imperial}.

\end{proof}

With the above results, we can prove the following two theorems.

\begin{teorema}
Let $(a,b)$ be an admissible pair. Then the operator $\mathcal P:L^1([0,1],\mu) \to L^1([0,1],\mu),$ admits  eigenvalue $\eta$ with respect eigenvalue $\lambda$ such that $\mu(\{\eta >0\})=1$ and $\|\eta\|_{L^1(M,\mu)}=1.$ In particular $Y_n^{a,b}$ fulfills hypothesis \ref{(K)}.\label{teorema612}
\end{teorema}
\begin{proof}
From Proposition \ref{vidaboa}, there exists an eigenmeasure $\nu(\d x) = h (x) \d x$ of $T^*$ with $h \in L^1([0,1],\mu).$ This implies that for every $f\in \mathcal C^0(M)$,
$$ \int_0^1 T(f)(x) h(x) \d x = \int_0^1 f(x) h (x) \d x.    $$
On the other hand, since $f g\in L^\infty([0,1])$ we obtain that
\begin{align*}
    \int_0^1 f(x) h (x) \d x &= \int_0^1 Tf (x) h(x) \d x = \int_0^1 \frac{\mathcal L (fg) (x)}{ \lambda g(x)} h(x) \d x= \int_0^1 f(x)  \frac{g(x)}{\lambda}\mathcal P\left( \frac{h}{g}\right) \d x.\\
\end{align*}
Finally, defining $\eta(x) = h(x)/g(x)$, it follows that $\mathcal P\eta = \lambda \eta.$ Since 
$$\eta (x) = \frac{1}{\lambda (b-a) x(1-x)} \int_{a x(1-x)}^{b x(1-x)\land 1} \eta (y) \d y, $$
we clearly have that $\eta \in \mathcal C^0( (0,1)).$ Moreover, it is easy to see that if there exists  $x_0\in (0,1)$ such that $ \eta (x_0) = 0,$ then $\eta( x) = 0$ $\mathrm{Leb}$-a.s. in $(0,1)$ which is a contradiction.
\end{proof}

\begin{teorema} \label{Thm:ABLogistic}
Let $(a,b)$ be an admissible pair. Consider $M=[0,1]$ and the Markov chain $Y_{n+1}^{(a,b)} = \omega_n Y_n^{(a,b)} (1-Y_n^{(a,b)})$ absorbed at $\partial = \R \setminus M,$ with $\{\omega_n\}_{n\in\mathbb N}$ an i.i.d sequence of random variable such that $\omega_n\sim \mathrm{Unif}([a,b])$ on $\mathbb R_{M}$ with absorption $\partial$. Then
\begin{enumerate}
    \item[${(i)}$] $Y_n^{a,b}$ admits a quasi-stationary measure $\mu_{a,b}$ with survival rate $\lambda_{a,b}$ such that $\supp(\mu) =[0,1]$ and $\mu \ll \mathrm{Leb}$, where  $\mathrm{Leb}$ denotes the Lebesgue measure on $[0,1].$
    \item[${(ii)}$] There exists $\eta^{a,b} \in L^1(M,\mu)$ such that $\mathcal P \eta^{a,b} = \lambda_{a,b} \eta^{a,b} $, $\|\eta^{a,b}\|_{L^1(M,\mu)} =1$ and $\eta^{a,b} >0$ $\mu_{a,b}$-almost surely.
    \item[${(iii)}$]
     For every $h\in L^\infty(M,\mathrm{Leb}),$
    $$ \lim_{n\to\infty} \mathbb E\left[ \frac{1}{n}\sum_{i=0}^{n-1} h\circ Y_n^{a,b} \mid \tau >n\right] = \int_M h(y) \eta^{a,b}(y) \mu_{a,b}(\d y)\ \text{for every }x\in(0,1).$$
    \item[${(iv)}$] For every $h\in L^\infty (M,\mu)$ 
    $$\lim_{n\to \infty} \mathbb E_x\left[ h \circ X_i\mid \tau >n\right]  = \int h(y) \mu(\d y)\ \mbox{for every $x\in (0,1)$}. $$
\end{enumerate}
\end{teorema}
\begin{proof} 

Note that Theorem \ref{teorema612} implies that $Y_n^{(a,b)}$ satisfies Hypothesis items $\mathrm{(H1a)}$ and $\mathrm{(H1b)}$ of Hypothesis \ref{(K)}, also items $\mathrm{(H1c)}$ and $\mathrm{(H1d)}$ of Hypothesis \ref{(K)} follows from Propositions \ref{uk1} and \ref{imperial}. 

Once again, from Propositions \ref{uk1} and \ref{imperial} we obtain that $Y_n^{a,b}$ satisfies Hypothesis \ref{(R)} defining $K_i := [1/i,1-1/i],$ for every $i\in\mathbb N$. Also, since the logistic map $4 x(1-x)$ is chaotic in $[0,1]$, and $f:\mathbb R\times[a,b]\to\mathbb R,$ $f(x,\omega)=\omega x(1-x)$ is a continuous function, we conclude that $m=1$ in Theorem \ref{teorema2}. Therefore, the conclusions of the theorem follow directly from Theorem \ref{teorema4}.
\end{proof}

\subsection{Analysis of the admissible pairs}
\label{logistic}

Fixing a pair $(a,b) \in (0,4)\times (4,\infty)$ it is relatively easy to check if $(a,b)$ is an admissible pair or not. However, it is complicated to solve inequalities \eqref{ineq1} and \eqref{ineq2} in terms of $(a,b)$. In this section, we prove that for every $(a,b)\in [1,4)\times (4,\infty)$ is an admissible pair.

We start showing that for every $(a,b)\in [1,2)\times (4,\infty),$ inequality \eqref{ineq1} is fulfilled.
\begin{proposicao}
For every $(a,b)\in [1,2)\times (4,\infty)$ we have that
$$0\leq \frac{1}{2}-\frac{1}{2} \sqrt{1-\frac{2}{b}\left(1-\sqrt{1-\frac{4x}{a}
   }\right)} \leq \frac{a}{4},\ \text{for every }x\in\left[\frac{4 a^2 - a^3}{16},\frac{a}{4}\right].$$\label{prop614}
\end{proposicao}
\begin{proof}
Note that $\frac{1}{2}-\frac{1}{2} \sqrt{1-\frac{2}{b}\left(1-\sqrt{1-\frac{4x}{a}}\right)}$ is an increasing function in $x.$ Therefore, for every $x\in [(4a^2 - a^3)/16,a/4]$ we obtain
\begin{align*}
       0\leq & \frac{1}{2}-\frac{1}{2} \sqrt{1-\frac{2}{b}\left(1-\sqrt{1-\frac{4x}{a}}\right)} \leq \frac{1}{2}-\frac{1}{2} \sqrt{1-\frac{2}{b}}\leq \frac{1}{4}\leq \frac{a}{4}.
\end{align*}   

\end{proof}

\begin{proposicao}
For every $b>4$, the following maps 
\begin{align*}
F^{a,b}_1(x):= \frac{2\left( \tanh^{-1}\left(\sqrt{1 - \frac{2}{b} + \frac{2}{b} \sqrt{1-\frac{4x}{b}}}\right) -  \tanh^{-1}\left(\sqrt{1 - \frac{2}{a} +\frac{2}{a}\sqrt{1-\frac{4 x}{b}}}\right)\right)   }{\sqrt{1-\frac{4x}{b}}} 
\end{align*}
and
\begin{align*}
  F^{a,b}_2(x):=  \frac{2 \tanh ^{-1}\left(\sqrt{\frac{-2 +2 \sqrt{1-\frac{4 x}{a}}+b}{b}}\right)+\log\left(\frac{a}{4-a}\right)}{\sqrt{1-\frac{4x}{a}}}
\end{align*}
are increasing in $x$ in the interval $[(4a^2 -a^3)/16,a/4].$\label{prop615}
\end{proposicao}
\begin{proof}
It is readily verified that $${x\mapsto  \tanh^{-1}\left(\sqrt{1 - \frac{2}{b} + \frac{2}{b} \sqrt{1-\frac{4x}{b}}}\right) -  \tanh^{-1}\left(\sqrt{1 - \frac{2}{a} +\frac{2}{a}\sqrt{1-\frac{4 x}{b}}}\right)  }\ \text{and }x \mapsto \sqrt{1-\frac{4x}{b}}$$ are, respectively, increasing and decreasing for $x\in [(4a^2 -a^3)/16,a/4],$ implying that $F^{a,b}_1(x)$ is increasing in $x$ in the interval $[(4a^2 -a^3)/16,a/4].$ 

In the following we prove that the $F^{a,b}_2$ is an increasing function in $[(4a^2 -a^3)/16,a/4].$ Through the change of coordinates $y =\sqrt{\frac{-2 +2 \sqrt{1-4x/a}+b}{b}},$ we obtain that to show that $F^{a,b}_2$ is an increasing function it is enough to show that
$$F^{a,b}_3(y) = \frac{\displaystyle  \log \left(\frac{1+y}{1-y}\right) + \log\left(\frac{a}{4-a}\right)}{b y^2-b+2}$$
is decreasing in $x$ in the interval $\left[\sqrt{(b-2)/b} ,\sqrt{(b-1)/b}\right]\supset \left[\sqrt{(-2+b)/b} ,\sqrt{(b-a)/b}\right].$ Since,
$$\frac{\d F^b_3}{\d y}(y) = \frac{2 \left(b \left(y^2-1\right) y \left(\log \left(\frac{a}{4-a}\right)+\log
   \left(\frac{1+y}{1-y}\right)\right)+b y^2-b+2\right)}{\left(1-y^2\right) \left(b
   \left(y^2-1\right)+2\right)^2},$$
it is enough to show that 
\begin{align*}
    b \left(y^2-1\right)\left( y \log
   \left(\frac{a}{4-a} \cdot \frac{1+y}{1-y}\right)+1\right)+2 &\leq  b \left(y^2-1\right) \left(y\log
   \left(  \frac{1+y}{3-3 y}\right)+ 1\right)+2 \\
   &\leq 0,\ \text{for every }y\in\left[\sqrt{\frac{b-2}{b}} , \sqrt{\frac{b-1}{b}}\right]. 
\end{align*} 
Observe that given $y\in \left[\sqrt{(-2+b)/b} ,\sqrt{(-1+b)/b}\right]\subset [\sqrt{2}/2,1],$ we obtain that
$$ b y (y^2-1) \log\left(\frac{1+ y}{3- 3y}\right) \leq 0\ \text{ and }  \frac{1+y}{3-3y} -1 \geq 0.$$

From \cite[Equation (2)]{log} it follows that $\log(1+x)\geq x/(1+x/2)$ for every $x\geq 0.$ Therefore
\begin{align}
    2+ b \left(y^2-1\right)+b y \left(y^2-1\right) \log \left(\frac{y+1}{3-3
   y}\right) &\leq 2+ b \left(y^2-1\right)+\frac{b y \left(\frac{y+1}{3-3 y}-1\right) \left(y^2-1\right)}{\frac{1}{2}
   \left(\frac{y+1}{3-3 y}-1\right)+1}\nonumber\\
   &=2-\frac{b (1-y^2) (4 y^2-3y+2)}{2-y}.\label{altralis1}
\end{align}
Using standard techniques, one can check that ${b>4}$ and $y \in [\sqrt{(-2+b)/b} $,$\sqrt{(-1+b/b}]$ then \eqref{altralis1} is less or equal to $0$, implying that $F^{a,b}_3$ is decreasing in the interval $\left[\sqrt{(-2+b)/b} ,\sqrt{(-1+b)/b}\right]$ for every $(a,b) \in [1,2]\times (4,\infty)$and therefore $F^{a,b}_2$ is increasing in $[(4a^2 -a^3)/16,a/4],$ for every $(a,b)\in [1,2]\times (4,\infty)$.

\end{proof}

Using the above proposition, we show that if $(a,b) \in [1,4)\times (4,\infty)$, then $(a,b)$ is an admissible pair for every.
\begin{teorema}
If $(a,b) \in[1,4)\times (4,\infty)$, then $(a,b)$ is an admissible pair.\label{T616}
\end{teorema}

\begin{proof}
From the definition of admissible pair, we just need to consider the case $(a,b) \in [1,2]\times (4,\infty).$ From Proposition \ref{prop614} we obtain that the pair $(1,b)$ satisfies equation \eqref{ineq1}.  

In the following, we show that the pair $(1,b)$ satisfies \eqref{ineq2}. Observe that $\eqref{ineq2}$ is equivalent of showing that $F^{a,b}_1(x) \leq F^{a,b}_2(x)$ for every $x\in[(4a^2 - a^3)/16,a/4]$,
where $F^{a,b}_1$ and $F^{a,b}_2$ are defined in Proposition \ref{prop615}. From Proposition \ref{prop615} it is enough to show that $$F_1^{a,b}\left(\frac{a}{4}\right)  \leq F_2^b\left(\frac{4a^2 -a^3}{16}\right) \text{ for every }(a,b)\in [1,4)\times(4,\infty).$$

We divide the proof into two steps.

\begin{step}{1}
\it We show that for every $b>4$, $(1,b)$ is an admissible pair.
\end{step}
Note that for every $b>4,$
\begin{align*}
    F_1^{1,b}\left(\frac{1}{4}\right) &=\frac{2 \left(\tanh ^{-1}\left(\sqrt{\frac{b+2
   \sqrt{1-\frac{1}{b}}-2}{b}}\right)-\tanh ^{-1}\left(\sqrt{2
   \sqrt{1-\frac{1}{b}}-1}\right)\right)}{\sqrt{1-\frac{1}{b}}}\\
   &\leq 4 \tanh ^{-1}\left(\frac{\sqrt{2 \sqrt{\frac{b-1}{b}}-1}-\sqrt{\frac{b+2
   \sqrt{\frac{b-1}{b}}-2}{b}}}{1- \sqrt{\left(2 \sqrt{\frac{b-1}{b}}-1\right)\left(\frac{b+2
   \sqrt{\frac{b-1}{b}}-2}{b}\right)}}\right)
\end{align*}
and
$$ F_2^{1,b}\left(\frac{3}{16}\right) = 4 \tanh ^{-1}\left(\frac{\sqrt{1-\frac{1}{b}}-\frac{1}{2}}{1-\frac{1}{2}
   \sqrt{1-\frac{1}{b}}}\right) = 4 \tanh^{-1}\left(\frac{3 \sqrt{b-1} \sqrt{b}-2}{3 b+1}\right).$$

Since the function $x\mapsto 4 \tanh^{-1}(x)$ is increasing, to finish the proof of the theorem, it is enough to show that
\begin{align}
    \frac{\sqrt{-1+2 \sqrt{\frac{b-1}{b}}}-\sqrt{\frac{b+2
   \sqrt{\frac{b-1}{b}}-2}{b}}}{1- \sqrt{\left(-1+2 \sqrt{\frac{b-1}{b}}\right)\left(\frac{b+2
   \sqrt{\frac{b-1}{b}}-2}{b}\right)}} \leq \frac{3 \sqrt{b-1} \sqrt{b}-2}{3 b+1}\ \text{for every }b>4. \label{ineq3}
\end{align}
Using standard methods, one can show that the above equation simplifies in showing that
$$p(b):=4239 b^6-23868 b^5+31482 b^4+8964
   b^3-40401 b^2+23424 b -4096 \geq 0\ \text{for every }b>4.$$
However, since for every $\delta>0$
$$p(4+\delta) = 4239 \delta ^6+77868 \delta ^5+571482 \delta ^4+2119716 \delta ^3+4091679 \delta ^2+3683256 \delta
   +998384>0, $$
we obtain that \eqref{ineq3} holds. This completes Step 1.

\begin{step}{2}
\it We show that $(a,b)$ is an admissible pair for every $(a,b)\in [1,2)\times(4,\infty).$
\end{step}
Fixed $b>4,$ observe that
$$(2-a) F_1^{a,b}\left(\frac{a}{4}\right) = 2\frac{2-a}{\sqrt{1-\frac{a}{b}}} \left(\tanh ^{-1}\left(\sqrt{\frac{2 \sqrt{1-\frac{a}{b}}+b-2}{b}}\right)-\tanh
   ^{-1}\left(\sqrt{\frac{2 \sqrt{1-\frac{a}{b}}+a-2}{a}}\right)\right)   $$
and 
$$ (2-a) F_2^{a,b}\left(\frac{4a^2 - a^4}{16}\right) = 2 \left(2 \tanh ^{-1}\left(\sqrt{\frac{b-a}{b}}\right)+\log
   \left(\frac{a}{4-a}\right)\right).$$

It is readily verified that
$$\frac{2-a}{\sqrt{1-\frac{a}{b}}}\ \text{and }\tanh ^{-1}\left(\sqrt{\frac{2 \sqrt{1-\frac{a}{b}}+b-2}{b}}\right)-\tanh
   ^{-1}\left(\sqrt{\frac{2 \sqrt{1-\frac{a}{b}}+a-2}{a}}\right) $$
are decreasing functions in $a\in [1,2),$ implying that $(2-a) F_1^{a,b}(a/4)$ is a decreasing function in $a\in[0,1]$ and 
$$ (a-2) F_2^{a,b} ((4a^2 - a^3)/16) = 2 \left(2 \tanh ^{-1}\left(\sqrt{\frac{b-a}{b}}\right)+\log
   \left(\frac{a}{4-a}\right)\right).$$
is an increasing function in $a \in [1,2)$.   From Step $1$ we obtain that for every $a\in [1,2),$
$$F_1^{a,b}\left(\frac{a}{4}\right)=\frac{(2-a)F_1^{a,b}\left(\frac{a}{4}\right)}{2-a} \leq \frac{ F_1^{a,b}(1/4)}{2-a}\leq \frac{F_2^{a,b}(3/16)}{2-a} \leq F_2^{a,b}\left(\frac{4a^2-a^3}{16}\right).$$
This completes the proof of the theorem.
\end{proof}

We finish the paper proving Theorem \ref{Thm:Logistic}.

\begin{proof}[Proof of Theorem \ref{Thm:Logistic}]

The theorem follows directly from Theorems \ref{Thm:ABLogistic} and \ref{T616}
\end{proof}

%  \includepdf[pages=-]{Prova.pdf}

\section*{Acknowledgments}

We are grateful to %Prof. Dr. 
Jochen Glück for many useful discussions and valuable comments regarding the theory of positive integral operators. MC’s research has been supported by an Imperial College President’s PhD scholarship. MC, VG and JL are also supported by the EPSRC Centre for Doctoral Training in Mathematics of Random Systems: Analysis, Modelling and Simulation (EP/S023925/1). JL gratefully acknowledges research support from IRCN, University of Tokyo and CAMB, Gulf University of Science and Technology, as well as from the London Mathematical Laboratory.

\bibliographystyle{plain} 
\bibliography{main}
\end{document}